\documentclass[12pt]{article}
\usepackage[english]{babel}
\usepackage[latin1]{inputenc}
\usepackage{amsmath}
\usepackage{amsfonts}
\usepackage{amssymb}
\usepackage{a4wide}

\def\Netoile{\mathbb{N} ^\ast}

\def\Rplusetoile{\R _+ ^\ast}

\newcommand{\N}{\mathbb{N}}

\newcommand{\R}{\mathbb{R}}

\newcommand{\Z}{\mathbb{Z}}

\newcommand{\al}{\alpha}

\newcommand{\congru}{\equiv}

\newcommand{\eps}{\varepsilon}

\newcommand{\findem}{\end{pf}}

\newcommand{\inclus}{\subset}

\newcommand{\moins}{\setminus}

\newcommand{\ro}{\varrho}

\newtheorem{Th}{Theorem}[section]
\newtheorem{Prop}[Th]{Proposition}
\newtheorem{Cor}[Th]{Corollary}
\newtheorem{Lemme}[Th]{Lemma}

\newtheorem{Defith}[Th]{Definition} 
\newtheorem{Example}[Th]{Example} 
\newtheorem{Remark}[Th]{Remark} 
 
\def\Dem{\hspace{-\parindent}\textsc{Proof }}
\def\Demdeuxpoints{\hspace{-\parindent}\textsc{Proof: }}

\newcommand{\specqcq}{{\cal D}} 
\newcommand{\spectre}{{\cal D}_0} 
\newcommand{\specBL}{{\cal D}'} 

\newcommand{\poubelleprovisoire}[1]{}

\newcommand{\inter}{\cap}

\newcommand{\miroir}[1]{\widetilde{#1}}

\newcommand{\eneq}{\end{equation}}
\newcommand{\begineq}{\begin{equation}}

\newcommand{\densi}{\delta}

\newcommand{\periw}{\pi}
\newcommand{\lgperiw}{d}
\newcommand{\unlgperiw}{\{1,\ldots,\lgperiw\}}

\newcommand{\qaa}{\alpha}
\newcommand{\pal}{\pi}
\newcommand{\imin}{i_0}
\newcommand{\fctpsi}{\psi}
\newcommand{\fctpsipr}{\psi '}
\newcommand{\lgr}[1]{\vert #1 \vert}
\newcommand{\alphb}{{\cal A}}

\newcommand{\borne}{B}
\newcommand{\difn}{\eps}
\newcommand{\quo}{\alpha}
\newcommand{\distamoins}{\delta^{-}}
\newcommand{\distaplus}{\delta^{+}}

\newcommand{\motdepsipr}{w_{\fctpsi'}}
\newcommand{\motdepsi}{w_{\fctpsi}}
\newcommand{\motdepsin}{w_{\fctpsi_n}}
\newcommand{\motvide}{\varepsilon}
\newcommand{\undisj}{\sqcup}

\newcommand{\jpri}{i}
\newcommand{\doublerightarrow}{\Rightarrow}

\newcommand{\app}{abundant palindromic prefixes} 
\newcommand{\appmajusc}{Abundant Palindromic Prefixes} 
\newcommand{\sig}{\sigma}
\newcommand{\caral}{c_{\alpha}}
\newcommand{\palstu}{\widehat{\pi}} 
\newcommand{\ome}{\omega}

\def\qed{\unskip\kern 8pt\penalty 500\raise 0pt\hbox{
\vrule\vbox to 6pt{\hrule width 4pt\vfill\hrule }\vrule }\par}

\title{Palindromic Prefixes and Episturmian Words}
\author{St\'ephane Fischler}
\date{\today}

\begin{document}

\maketitle

{\bf Abstract:} Let $w$ be an infinite word on an alphabet $\alphb$. We
denote
by $(n_i)_{i \geq 1}$ the increasing sequence
(assumed to be infinite)
 of all lengths of palindromic prefixes of $w$. In this text, we give an explicit construction of
all words $w$ such that $n_{i+1} \leq 2 n_i + 1$ for all $i$, and study these 
words. Special examples include characteristic Sturmian words, and more
generally standard episturmian words. As an application, we 
study the values taken by
the quantity $\limsup n_{i+1}/n_i$, and prove that it is minimal (among
all non-periodic words) for the Fibonacci word.
 
\bigskip

\section{Introduction}

The purpose of this text is to study infinite words 
(on an arbitrary, not necessarily finite, alphabet $\alphb$)
which have ``sufficiently many'' palindromic prefixes.
The motivation comes from diophantine approximation (see below), though this question is also related
to physics, namely to the spectral theory of discrete one-dimensional Schr\"odinger 
operators. Words with many palindromic factors can be used in this setting 
\cite{HofKnillSimon}, corresponding to the combinatorial notion of ``palindrome complexity'' (see for
instance \cite{ABCD}). On the other hand, replacing ``whole-line methods'' by ``half-line methods'' 
in connection with this problem 
leads \cite{DamanikGhezRaymond} to the use of words with many palindromic prefixes, like the ones
studied below.

\bigskip

In precise terms, given an infinite word $w$, we shall denote (in this Introduction)
by $(n_i)_{i \geq 1}$ the increasing sequence of all lengths of palindromic prefixes of $w$, with 
$n_1 = 0$ corresponding to the empty prefix. The words studied here always have infinitely many 
palindromic prefixes, so we assume the sequence $(n_i)_{i \geq 1}$ to be infinite.

A trivial example of such a word $w$ is any periodic word with a palindromic period. A more interesting example 
is the Fibonacci word $w = babbababbabba \ldots$ 
on the two-letter alphabet $\{a,b\}$,
 for which the sequence $(n_i) = (0,1,3,6,11,\ldots)$ is given by $n_i = F_{i+1}-2$ (where 
$F_i$ is the $i$-th Fibonacci number); this follows from \cite{DeLuca} (Theorem 5).
  More generally, any characteristic Sturmian word satisfies
$n_{i+1} \leq 2 n_i + 1$ for any $i$, and denoting by $[0,s_1,s_2,\ldots,]$ the continued
fraction expansion of its slope we have (see \S \ref{subsecstu}): 
\begin{equation} \label{eqstu}
\limsup \frac{n_{i+1}}{n_i} = \limsup [1,1,s_k,s_{k-1},\ldots,s_1]. 
\eneq
In particular, if $w$ is the Fibonacci word then 
$\limsup n_{i+1}/n_i$ is the golden ratio $\gamma = (1+\sqrt5)/2$.

\medskip

A generalization of characteristic Sturmian words to an arbitrary alphabet has been given by Droubay, Justin and Pirillo \cite{DJP}:
these are standard episturmian words. They also satisfy $n_{i+1} \leq 2 n_i + 1$ for any $i$,
but there is no easy equation like \eqref{eqstu} to compute $\limsup n_{i+1}/n_i$. 

\medskip

In this text, we study the words $w$ 
with abundant palindromic prefixes
in the following sense:

\begin{Defith} An infinite word $w$ is said to have {\em 
abundant palindromic prefixes} if the sequence $(n_i)_{i \geq 1}$ of all lengths of
its palindromic prefixes is infinite and satisfies 
$n_{i+1} \leq 2 n_i + 1$ for any $i \geq 1$.
\end{Defith}

A completely explicit construction of all
words with abundant palindromic prefixes is given, which generalizes one of the 
constructions \cite{JP} 
of standard episturmian words.
 This is a strict generalization, i.e., there are words 
with \app\ which 
are not standard episturmian. Moreover, our results extend to words such that
$n_{i+1} \leq 2 n_i + 1$ for any sufficiently large integer $i$; in particular, a general
construction of all such words is given.

\bigskip

For any word $w$, we let 
$$\densi(w) = \limsup \frac{n_{i+1}}{n_i}$$
if $w$ admits infinitely many palindromic prefixes, and $\densi(w) = \infty$ otherwise.
Then $1/\densi(w)$ measures the ``density'' of palindromic prefixes in $w$. We let $\specqcq$ be the set
of real numbers that can be written $\densi(w)$ for some word $w$
(on a suitable alphabet).
Moreover, we let $\spectre$ be the set of all numbers $\densi(w)$ obtained
from words $w$ with \app.
 The inclusion
$\spectre \inclus \specqcq \inter [1,2]$ trivially holds, and it is not difficult to prove
(see \S \ref{subsecscarse}) that $(2, +\infty] \inclus \specqcq$. Denoting by $\undisj$ the union of
two disjoint sets, the following result holds.

\begin{Th} \label{theoqcq}
We have $\specqcq = \spectre \undisj (2,+\infty]$.
\end{Th}

Actually, for any word $w$ such that $\densi(w) < 2$, there is a 
word $w'$ with \app\ such that the palindromic prefixes of $w$ satisfy the same
recurrence relation as those of $w'$ (see
Proposition \ref{propoutilqcq}) and, therefore, $\densi(w) = \densi(w')$. 

\bigskip

The easiest examples of words with \app\ are periodic words
(with a palindromic period) and 
characteristic Sturmian words (for which $\densi(w)$ can be computed
thanks to Equation \eqref{eqstu}). Denote 
by $\specBL$ the set of numbers $\densi(w)$, for these words $w$. Obviously we have 
$\specBL \inclus \spectre$, and the following theorem shows that this inclusion is an equality if we restrict
to words with ``sufficiently many'' palindromic prefixes:

\begin{Th} \label{theosqrttrois}
We have $\specBL \inter [1,\sqrt3] = \spectre \inter [1,\sqrt3] = \specqcq \inter [1,\sqrt3] $.
\end{Th}

For a periodic word $w$ with a palindromic period, 
we have trivially $\densi(w) = 1$. For a characteristic Sturmian word $w$ with slope
 $[0,s_1,s_2,\ldots,]$, Equation \eqref{eqstu} allows one to compute $\densi(w)$.
From this it is easy to deduce that
the characteristic Sturmian word $w$ with minimal value of $\densi(w)$ is 
the Fibonacci word. This shows that $1$ and the golden ratio $\gamma = (1+\sqrt5)/2$ are the two
smallest elements in $\specBL$. Cassaigne studied 
(\cite{Cassaigne}, Corollary 1 and Theorem 2)
the next elements, and his result (together
with Theorem \ref{theosqrttrois}) yields:

\begin{Th} The smallest elements in $\spectre$ (resp. in $ \specqcq$)
make up an increasing sequence $(\sigma_n)_{n \geq 0}$ with 
$\sigma_0 = 1$ and $\sigma_1 = \gamma$, converging to the smallest accumulation point
$\sigma_{\infty}$ of $\spectre$ (resp. of $ \specqcq$).
\end{Th}

In more precise terms, this statement means
that $\spectre \inter [1,\sigma_{\infty}) = \{\sigma_n, n \geq 0\}$.
Moreover all $\sigma_n$, and $\sigma_{\infty} = 1.721\ldots$, are given 
in an explicit way in terms of their continued
fraction expansion. For instance, writing $\overline m$ for the periodic repetition
$mmm\ldots = m^\ome$ of a finite sequence $m$, we have:
\begin{eqnarray*}
\sigma_2 &= 1+\sqrt{2}/2 &= 1.707\ldots = [1,1,\overline{2}]\\
\mbox{ and } \sigma_3 &= (2+\sqrt{10})/3 &= 1.720 \ldots = [1,1,\overline{2,1,1}].\\
\end{eqnarray*}

As a corollary, we see that the Fibonacci word has maximal ``palindromic prefix density'' among non-periodic words:

\begin{Cor} \label{corfibo}
Let $w$ be an infinite word with $\densi(w) < \gamma$. Then $w$ is periodic.
\end{Cor}

For a characteristic Sturmian word $w$ with slope
 $[0,s_1,s_2,\ldots,]$, Morse 
and Hedlund have computed \cite{MorseHedlund} the recurrence function of $w$.
This gives (see Corollary 1 of \cite{Cassaigne}) a formula for the 
 recurrence quotient $\ro(w)$ of $w$, namely
 $\ro(w) = 2 + \limsup [s_k,s_{k-1}, \ldots, s_1]$.
Therefore Equation \eqref{eqstu} gives in this case:
$$\densi(w) = \frac{2 \ro(w) - 3}{\ro(w) - 1}, $$
hence (as above) the Fibonacci word has minimal recurrence quotient (equal to 
$(5+\sqrt5)/2$) among all
characteristic Sturmian words. Rauzy has conjectured \cite{RauzyBordeaux} that it 
has minimal recurrence quotient among all non-periodic words. Corollary \ref{corfibo} is an analogue of
this conjecture.  

\bigskip

The motivation for this text comes from diophantine approximation. Actually $\spectre
\moins \{1\}$ is
equal \cite{SFOttdio} to the
set denoted by ${\cal S}_0 \inter [1,2]$ in \cite{CRASasim}, defined in terms of an exponent that
measures the simultaneous approximation to a real number and its square by rational numbers
with the same denominator. In particular, Theorem 2.1 in \cite{CRASasim} follows from this equality
and Theorem \ref{theosqrttrois} stated above.

This connection between palindromic prefixes and diophantine approximation is due to 
Roy \cite{RoyCRAS}. It allows
one  
to get a purely number-theoretical proof of Corollary \ref{corfibo} stated above, by applying Davenport-Schmidt's theorem 
\cite{DS} on simultaneous approximation to $\xi$ and $\xi^2$ to the real number $\xi$
obtained (as in \cite{RoyCRAS}) from an infinite word $w$.

\bigskip

The structure of this text is as follows. 
We first explain the notation (\S \ref{subsecnotation}), and prove that
for any $\al > 2$ there is a word $w$ such that $\densi(w) = \al$ (\S \ref{subsecscarse}).
This explains why the rest of the text is devoted only to words $w$ such that 
$\densi(w) \leq 2$.
 
Then we recall how
characteristic Sturmian words (\S \ref{subsecstu}) and
standard episturmian words (\S \ref{subsecconstrepistu}) are constructed, with a special emphasis
on their palindromic prefixes. In Section \ref{secapp}, we construct 
all words with \app\ (\S \ref{subsecconstruapp}). To study these words, the
key definition is the one of reduced functions, which allows us to state
(\S \ref{subsecreduced}) the main results on words with \app . Moreover, we
explain (\S \ref{subsecdensipsi}) how to compute $\densi(w)$ for such a word $w$, using
the associated reduced function $\fctpsi$. 
The proof of the results stated in Section \ref{secapp} is given in Section \ref{secpreuve},
using general lemmas (\S \ref{subseclemmesgeneraux} and \ref{subsecperiodq}) that might be of
independent interest.

Next we briefly explain how to generalize the results of Section \ref{secapp}
to words that satisfy $n_{i+1} \leq 2n_i + 1$ for any sufficiently
large $i$ (\S \ref{subsecasyabondant}). This allows to prove (in \S \ref{subsecdemtheoqcq})
Theorem \ref{theoqcq} stated above.

Theorem \ref{theosqrttrois} is proved in \S \ref{subsecpreuvetheosqrttrois}, and the set
$\spectre$ (resp. $\specBL$) is studied near $\sqrt3$ in \S \ref{subsecracinetroisnonepistu}
(resp. in \S \ref{subsecspectresqrt3}); this implies that Theorem \ref{theosqrttrois} is 
optimal. We also define $\alphb$-strict words with \app\
in \S \ref{subsecinitial}, and prove that any $\alphb$-strict standard
 episturmian word $w$ such that
$\densi(w) < \sqrt3$ is either periodic or characteristic Sturmian.

Section \ref{secquestions} contains questions and open problems about words
with \app . At last,
 Section \ref{secappendice} is an appendix 
 devoted to the proof of two technical results: Proposition \ref{propindep} (stated
 in \S \ref{subsecdemtheoqcq}) and Lemma \ref{lemdecalcul} 
 (stated in \S \ref{subsecpreuvetheosqrttrois}). These statements concern 
asymptotic properties of the sequence $(n_i)$ associated with
 a word $w$ such that $\densi(w) < 2$. They are also
useful for the diophantine analogue \cite{SFOttdio} of this text.
 
\bigskip

\hspace{-\parindent}{\bf Acknowledgements:} I am very thankful to Jean-Paul Allouche
and Boris Adamczewski for their help and for pointing out to me crucial references in 
combinatorics. I would like also to thank 
Damien Roy and Michel Waldschmidt for their support in the number-theoretic 
counterpart of this paper, and Jacques Justin
for many useful remarks. 
At last, I am indebted to Jimena Sivak for help in the redaction.

\section{Notation and a Peculiar Construction} \label{secdebut}

\subsection{Notation} \label{subsecnotation}

Throughout the text, we consider a (finite or infinite) alphabet
$\alphb$, which we assume to be disjoint from 
$\Netoile = \{1,2,3,\ldots\}$. Of course, this is not a serious
restriction; it allows us to consider $\alphb \undisj \Netoile$ as
a disjoint union.

We denote by $\lgr{u}$ the length of a finite word $u$, that is the number
of letters in $u$, and by $\motvide$ the empty word (which
has length zero).
Given a finite word $u = u_1\ldots u_p$ with $u_i \in \alphb$ for any
$i \in \{1, \ldots, p\}$, we denote by $\miroir{u}$ its mirror image
$u_p \ldots u_1$, in such a way that $u$ is a palindrome if, and only if, $u = \miroir{u}$.
We set  $\miroir{\motvide} = \motvide$, so that  $\motvide$ is considered a palindrome.
We say that a word $u' = u'_1\ldots u'_{p'}$ is a prefix of $u$ if 
$p' \leq p$ and 
$u_j = u'_j$ for any $j \leq p'$, that is if there is a word $u''$ 
such that $u = u'u''$. We extend this definition to the case where
$u$ is an infinite\footnote{In this text, we consider only right
infinite words. In particular, all palindromes are assumed to be
finite.} word $u_1u_2u_3\ldots$.
In particular, $\motvide$ is a palindromic prefix of any 
(finite or infinite) word.

In the same way, a word $u''$ is a suffix of $u$ if, and only if,
there is a finite word $u'$ such that $u=u'u''$. If this happens then
either both $u$ and $u''$ are finite, or both $u$ and $u''$
are infinite.

\medskip

If $w$ and $w'$ are finite 
words such that $w'$ is a prefix of $w$, we denote by ${w'}^{-1}w$ the word $w''$ 
such that $w = w' w''$. In the same way, if $w = w'w''$, we write $w' = w {w''}^{-1}$.
An important special case is the following: if $w$ and $w'$ are palindromes and $w'$ is a prefix of $w$, then
$w'$ is also a suffix of $w$ and $w {w'}^{-1} w$ is again a palindrome (of which $w$ is a prefix). 
In this situation, if $w = w' w''$ then we have $w {w'}^{-1} w = w' {w''}^2$
(see Lemma \ref{lemfacile} below).

\begin{Remark} Let $w$ be a word on the (finite or infinite) alphabet $\alphb$, such that $n_{i+1}
\leq 2 n_i$ for any $i$ sufficiently large (with the sequence $(n_i)_{i \geq 1}$ defined in the
Introduction). Then only finitely many letters of $\alphb$ occur in $w$; this follows from Proposition \ref{propoutilqcq}
proved below. Therefore the interesting case, throughout this paper, is when $\alphb$ is finite.
\end{Remark}

\subsection{Words with Scarce Palindromic Prefixes} \label{subsecscarse}

In this Section, we prove that $(2, +\infty] \inclus \specqcq$. This result explains why all words $w$ studied
in the rest of this text are such that $\densi(w) \leq 2$.

Obviously there are words $w$ with only a finite number of palindromic prefixes; they
satisfy $\densi(w) = \infty$ hence $\infty \in \specqcq$. Now
let $\al$ be a real number greater than 2, and choose $\eps > 0$ such that $2 +\eps < \al$.
Denote by $(p_k)_{k \geq 0}$ a sequence of positive integers such that 
$\frac{p_k}{10^k}$ tends to $\al$, with $\frac{p_k}{10^k} > 2$ for any $k$. We define an 
increasing sequence $(n_i)_{i \geq 1}$ in the following way. We let $n_1 = 0$, $n_2=1$ and
if $i \geq 2$ is even we let $v_{i+1}$ be the maximal integer such that there exists a
multiple of $10^{v_{i+1}}$, denoted by $n_{i+1}$, with
$2n_i < n_{i+1} < (2+\eps)n_i + 1$. If $i \geq 3$ is odd, we let
$n_{i+1} = p_{v_i} \frac{n_i}{10^{v_i}}$. With this definition, 
we have $n_{i+1} \geq 2n_i + 1$ for any $ i \geq 1$, and $v_i$ (which is defined only when $i$ is odd)
tends to infinity as $i$ tends to infinity. This implies $\limsup \frac{n_{i+1}}{n_i} = \al$.

Now let us construct a word $w$ such that $(n_i)_{i\geq 1}$ is exactly the sequence of all
lengths of palindromic prefixes of $w$. We consider an alphabet 
$\alphb = \{\delta_k, k \in \N\}$ with $\delta_i \neq \delta_j$ when $i \neq j$.
We define finite palindromes $\pal_i$, of length $n_i$, by
$\pal_1 = \motvide $ and, for $i \geq 1$:
$$\left\{
\begin{array}{l}
\pal_{i+1} = \pal_i \delta_i \delta_0 ^{n_{i+1}-2n_i-2} \delta_i \pal_i 
\mbox{ if } n_{i+1} \geq 2n_i + 2\\
\pal_{i+1} = \pal_i \delta_i \pal_i 
\mbox{ if } n_{i+1} = 2n_i + 1. 
\end{array}
\right.$$
Then for any $i \geq 1$, $\pal_i$ is a palindrome written on the alphabet
$\{\delta_0,\ldots,\delta_{i-1}\}$. It is also a prefix of $\pal_{i+1}$, and all
palindromic prefixes of $\pal_{i+1}$ (except $\pal_{i+1}$ itself) are prefixes
of $\pal_i$. The infinite word $w$ defined as the limit of $\pal_i$ as $i$ tends to infinity
satisfies the required property: its palindromic prefixes are exactly the $\pal_i$'s, with
$n_i = \lgr{\pal_i}$. Therefore $\al = \limsup \frac{n_{i+1}}{n_i} = \densi(w) \in \specqcq$.

This proves the desired result, namely $(2, +\infty] \inclus \specqcq$.

\begin{Remark} It is possible to adapt this construction to any fixed finite alphabet containing at 
least two letters. This proves that for any (finite or infinite) alphabet $\alphb$ with at least two letters,
and for any $\al \in (2, +\infty]$, there exists a word $w$ on the alphabet $\alphb$ such that 
$\densi(w) = \al$.
\end{Remark}

\section{Sturmian and Episturmian Words} \label{secstuepistu}

In this Section, we recall how to construct characteristic, or standard,  Sturmian
(\S \ref{subsecstu}) and standard episturmian (\S \ref{subsecconstrepistu})
words, with a focus on the properties of their palindromic prefixes.

\subsection{Characteristic Sturmian Words} \label{subsecstu}

In this Section, we recall a construction of characteristic Sturmian words (see \cite{Lothaire}, Chapter~2) and properties of their palindromic prefixes.

We consider the two-letter alphabet $\alphb = \{a,b\}$. 
Let $s_1$, $s_2$, \ldots, be an infinite sequence of positive integers.
Define $\sig_0 = a$, $\sig_1 = a^{s_1 - 1}b$ and, by induction,
$\sig_n = \sig_{n-1} ^{s_n} \sig_{n-2}$ for any $n \geq 2$. In the terminology
of \cite{Lothaire} (page~75), $(\sig_n)$ is the standard sequence associated with
$(s_1-1, s_2, s_3, \ldots)$. For any $n \geq 1$, $\sig_n$ is a prefix of $\sig_{n+1}$; therefore the
words $\sig_n$ tend to an infinite word $\caral$, called the {\em characteristic Sturmian word}
with slope $\alpha = [0,s_1,s_2,\ldots]$.

\bigskip

For $n \geq 2$ and $1 \leq p \leq s_n$, the word $\sig_{n-1}^p \sig_{n-2}$
is a prefix of $\sig_n = \sig_{n-1} ^{s_n} \sig_{n-2}$ (since $\sig_{n-2}$ is
a prefix of $\sig_{n-1}$), hence of $\caral$. Moreover it ends with $ba$ if $n$ is even, and 
with $ab$ if $n$ is odd. As it is a standard word, there exists a  palindrome $\palstu_{n,p}$
such that
\begin{equation} \label{eqdefpalstu}
\left\{
\begin{array}{l}
\sig_{n-1}^p \sig_{n-2} = \palstu_{n,p} ba \mbox { if $n$ is even}\\
\sig_{n-1}^p \sig_{n-2} = \palstu_{n,p} ab \mbox { if $n$ is odd.}
\end{array}
\right.
\end{equation}
Actually  $\palstu_{n,p}$ is even a central word, so it can be written $\pi a b \pi'$ for some palindromes $\pi$, $\pi'$ (see \cite{LucaMignosi}). However, 
in what follows, we shall use only the fact that
the words $\palstu_{n,p}$ defined in this way are palindromes. This fact
can be proved directly (see for instance \cite{BL}, Lemma 5.3). 

\bigskip

We shall now define a sequence $(\pal_i)_{i \geq 1}$ of palindromic prefixes of $\caral$. First, for any $k \geq 1$
we let 
$t_k = s_1 + \ldots + s_k$. Now observe that for 
any $i \geq s_1$ there is exactly
one pair $(n,p)$ with $n \geq 2$ and $1 \leq p \leq s_n$ such that
$i = t_{n-1} + p - 1$. Therefore the equality
$$\palstu_{n,p} = \pal_{t_{n-1} + p - 1} \mbox{ for $n \geq 2$ and $1 \leq p \leq s_n$}$$
defines $\pal_i$ in a unique way for $i \geq s_1$
(and $\pal_{t_k - 1} = \palstu_{k,s_k}$ is obtained from $\sig_k$ by removing the last two letters). If $s_1 \geq 2$, we let
$\pal_i = a^{i-1}$ for any $i \in \{1, \ldots, s_1-1\}$. Then $\pal_i$ is defined for any
$i \geq 1$; we have $\pal_1 = \motvide$ and each $\pal_i$ is a prefix of $\pal_{i+1}$. Moreover all
$\pal_i$'s are palindromic prefixes of $\caral$.

\bigskip

Actually the $\pal_i$'s are the only palindromic prefixes of $\caral$. This follows from
de Luca's result (\cite{DeLuca}, Theorem 5; see 
also \cite{DJP}, \S 3) that $\pal_{i+1}$ is the right palindromic closure 
of $\pal_i \delta_i$, where $\delta_i \in \alphb$ is the letter in $\pal_{i+1}$ that comes
right after $\pal_i$ (see \S \ref{subsecconstrepistu} below). Another proof of this result can
be obtained by applying Theorem \ref{aucunpalmanque} proved in this text (see Example \ref{exstuepistugal}).

Since the $\pal_i$'s are exactly the palindromic prefixes of $\caral$, we have the equality
$\densi(\caral) = \limsup \lgr{\pal_{i+1}}/\lgr{\pal_i}$. 
It is not difficult to deduce Equation \eqref{eqstu}
 from this (see \cite{BL}, Proposition 7.1).

\bigskip

Let $k \geq 3$. It is not difficult to prove the relation
\begin{equation} \label{equnintermstu}
\pal_{t_k + \ell} = \sig_k ^{\ell+1} \pal_{t_{k-1}-1} \mbox{ for any } 
\ell \in \{0,\ldots, s_{k+1}\}
\end{equation}
using (for the case $\ell = s_{k+1}$) the identity 
$\sig_{k-1} \pal_{t_k - 1} = \sig_k \pal_{t_{k-1}-1}$ (see for instance
\cite{BL}, Lemma 5.1). From Equation \eqref{equnintermstu} immediately follows
\begin{equation} \label{eqintermedstu}
\left\{
\begin{array}{l}
\pal_{t_k + 1} = \pal_{t_k} \pal_{t_{k-1}-1} ^{-1} \pal_{t_k} \\
\pal_{t_k + \ell + 1} = \pal_{t_k + \ell} \, 
			\pal_{t_k + \ell-1}^{-1} \, \pal_{t_k + \ell} 
					\mbox{ for any } 
					\ell \in \{1,\ldots, s_{k+1}-1\}. 
\end{array}
\right.
\end{equation}
We are going now to define a map $\fctpsi: \Netoile \rightarrow \Netoile
\undisj \alphb$ in such a way that,
 for any $i \geq 1$:
\begin{equation} \label{eqrecstu}
\left\{
\begin{array}{l}
\pal_{i+1} = \pal_i \pal_{\fctpsi(i)}^{-1} \pal_i \mbox{ if } 
			\fctpsi(i) \in \Netoile , \\
\pal_{i+1} = \pal_i \fctpsi(i) \pal_i \mbox{ if } 
				\fctpsi(i) \in \alphb. 
\end{array}
\right.
\end{equation}
The possibility to define inductively, in this way, the palindromic prefixes
of $\caral$ using $\fctpsi$ will be the crucial point in the
construction of Section \ref{secapp}.

For $k \geq 3$ we let $\fctpsi(t_k) = t_{k-1}-1$, and if $i > t_3$
is not among the $t_k$'s we let $\fctpsi(i) = i-1$. Then
Equation \eqref{eqintermedstu} shows that \eqref{eqrecstu}
holds for any $i \geq t_3$. To define the values $\fctpsi(i)$ for
$1 \leq i < t_3$, we distinguish between two cases.

First, let us assume $s_1 = 1$. Then $\pal_\ell = b^{\ell-1}$ for
$1 \leq \ell \leq t_2$ and $\pal_{t_2 + \ell} = (b^{s_2}a)^{\ell}b^{s_2}$
for any $0 \leq \ell \leq s_3$. We let $\fctpsi(1) = b$,
$\fctpsi(t_2) = a$ and $\fctpsi(i) = i-1$
for $i \in \{2, \ldots, t_3 - 1\} \moins \{t_2\}$. Then 
Equation \eqref{eqrecstu}
holds for any $i \geq 1$.

Now let us assume $s_1 \geq 2$. Then $\pal_\ell = a^{\ell-1}$ for
$1 \leq \ell \leq t_1$ and $\pal_{t_1 + \ell} = (a^{s_1-1}b)^{\ell}a^{s_1-1}$
for any $0 \leq \ell \leq s_2$. 
Moreover Equation \eqref{equnintermstu} holds also for $k=2$.
We let $\fctpsi(1) = a$,
$\fctpsi(t_1) = b$, $\fctpsi(t_2) = t_1 - 1$ and $\fctpsi(i) = i-1$
for $i \in \{2, \ldots, t_3 - 1\} \moins \{t_1,t_2\}$. Then 
Equation \eqref{eqrecstu}
holds for any $i \geq 1$.

\subsection{Standard Episturmian Words} \label{subsecconstrepistu}

Denote by $w^{(+)}$ the (right)
palindromic closure of a finite word $w$, that is the shortest palindrome of which
$w$ is a prefix.
Let $\Delta = \delta_1 \delta_2 \ldots$ be an infinite word on
an alphabet $\alphb$. Droubay, Justin and Pirillo gave \cite{DJP} the following definition (see \cite{JP}, Corollary 2.2):

\begin{Defith} \label{defiepistu}
The standard episturmian word
with directive word $\Delta$ is the limit of the sequence $(\pal_i)_{i \geq 1}$ defined by
 $\pal_1 = \motvide$ and $\pal_{i+1} = (\pal_i \delta_i)^{(+)}$
for $i \geq 1$. 
\end{Defith}
 
 The important point here (which will be generalized
in \S \ref{subsecconstruapp}) is that a standard episturmian word can be 
constructed as a limit of an infinite sequence
of its palindromic prefixes.

\bigskip 
 
Given $\Delta$, define a function $\fctpsi: \Netoile \rightarrow \Netoile \undisj \alphb$ as follows.
For $n \geq 1$, let $\fctpsi(n) = \delta_n$ if the letter $\delta_n$ occurs for the first time
in $\Delta$ at the $n$-th position. Otherwise, let $\fctpsi(n) = n'$ where $n'$ is the greatest integer
such that $1 \leq n' \leq n-1$ and $\delta_{n'} = \delta_n$. Then for any $i \geq 1$ 
we have (\cite{JP}, p. 287): 
$$\pal_{i+1} = \pal_i \pal_{\fctpsi(i)}^{-1} \pal_i \mbox{ if } \fctpsi(i) \in \Netoile$$
and 
$$\pal_{i+1} = \pal_i \fctpsi(i) \pal_i \mbox{ if } \fctpsi(i) \in \alphb. $$
The crucial remark in what follows is that these equalities could have been taken as a definition
of the sequence $(\pal_i)$, and therefore of standard episturmian words.

\begin{Example} 
Let $s_1,s_2,\ldots$ be a sequence of positive integers, and
 $\alphb = \{a,b\}$ be a two-letter alphabet.
The standard episturmian word 
 with directive word $\Delta = a^{s_1-1}b^{s_2}a^{s_3}b^{s_4}\ldots$
is the characteristic Sturmian word with slope $[0,s_1,s_2,\ldots]$.
This follows from \S \ref{subsecstu} (see also \cite{DeLuca}, proof
of Theorem 5).
\end{Example}

\begin{Example} Let $\alphb = \{a,b,c\}$ and $\Delta = (abc)^\ome =abc abc abc \ldots$. Then the 
standard episturmian word $w$ with directive word $\Delta$ is 
(\cite{JP}, Example 2.1)
the Tribonacci (or Rauzy \cite{RauzySMF}) word 
(that is,
the fixed point $abacabaabacabab\ldots$
of the morphism defined by 
$a \mapsto ab$, $b \mapsto ac$ and $c \mapsto a$). The corresponding
function $\fctpsi$ is given by $\fctpsi(n) = n-3$ for $n \geq 4$, and
$\fctpsi(n) = \delta_n$ for $1 \leq n \leq 3$.
\end{Example}

\section{Words with \appmajusc} \label{secapp}

In this Section, we give a general construction (\S \ref{subsecconstruapp}) of all words with
abundant palindromic prefixes, using functions $\fctpsi$. Then we define (\S \ref{subsecreduced})
{\em reduced} functions $\fctpsi$; this definition allows us to state the main results 
about words with \app , namely Theorems \ref{aucunpalmanque}
and \ref{interetreduit}. At last, we explain in \S \ref{subsecdensipsi}
how to compute $\densi(w)$ (for a word $w$ with \app ) using the associated 
reduced function $\fctpsi$.

\subsection{A General Construction} \label{subsecconstruapp}

Let 
$\fctpsi: \Netoile \rightarrow \Netoile \undisj \alphb$ be any 
map such that, for each $n \geq 1$:
$$\mbox{ either } \fctpsi (n) \in \alphb \mbox{ or } 1 \leq \fctpsi(n) \leq n-1. $$

Define $\pal_1 = \motvide$ and, for $i \geq 1$: 
$$\pal_{i+1} = \pal_i \pal_{\fctpsi(i)}^{-1} \pal_i \mbox{ if } \fctpsi(i) \in \Netoile$$
and 
$$\pal_{i+1} = \pal_i \fctpsi(i) \pal_i \mbox{ if } \fctpsi(i) \in \alphb. $$
It is not difficult to prove by induction that all $\pal_i$'s are palindromes, and that $\pal_i$ is a prefix of $\pal_{i+1}$ 
 (for instance, if $\fctpsi(i) \in \Netoile$, writing
$\pal_i = \pal_{\fctpsi(i)} b_i = \miroir{b_i} \pal_{\fctpsi(i)}$ yields
$\pal_{i+1} = \miroir{b_i} \pal_{\fctpsi(i)} b_i = \pal_{\fctpsi(i)} b_i ^2$; the easy 
Lemma \ref{lemfacile} stated below can also be used).
However, in general there is no letter $\delta_i \in \alphb$ such that
 $\pal_{i+1}$ be 
the palindromic closure of $\pal_i \delta_i$. 

\begin{Defith} \label{defigalepisu}
We
call {\em word with \app} associated with $\fctpsi$,
and denote by $\motdepsi$, the limit of the sequence
$(\pal_i)$.
\end{Defith}

This definition is consistent with the one given in the Introduction
since the following result holds (it is proved in Section \ref{secpreuve}
as a consequence of Theorem \ref{interetreduit} stated below):

\begin{Th} \label{thcarac} 
Let $w$ be an infinite word, and $(n_i)_{i \geq 1}$ 
be the increasing sequence (assumed to be infinite) of the lengths of
its palindromic prefixes (with $n_1 = 0$). Then the following statements are 
equivalent:
\begin{enumerate}
\item[$(i)$] We have $n_{i+1} \leq 2 n_i + 1$ for any $i \geq 1$ (i.e., $w$ has
\app ).
\item[$(ii)$] For some function $\fctpsi$, we have $w = \motdepsi$ (i.e., $w$ is the
word with \app\ associated with $\fctpsi$).
\end{enumerate}
\end{Th}

Let us study in more details the word with \app\ associated with 
a map $\fctpsi$. First, let us consider the letter $\delta_i$ in $\pal_{i+1}$ 
that comes right after $\pal_i$. This is the first letter of $\pal_i^{-1}\pal_{i+1}$,
the one such that $\pal_i \delta_i$ is a prefix of $\pal_{i+1}$. We have $\delta_i = \fctpsi(i)$
if $\fctpsi(i) \in \alphb$, and $\delta_i = \delta_{\fctpsi(i)}$ otherwise. This explains 
the following definition:

\begin{Defith} \label{defidirecfctpsi}
We call {\em word of first letters} associated with 
 $\fctpsi$ the word $\Delta = \delta_1 \delta_2 \ldots$ defined 
(for each $n \geq 1$) by
 $\delta_n = \fctpsi(n)$ if $\fctpsi(n) \in \alphb$, and
 $\delta_n = \delta_{\fctpsi(n)}$ otherwise.
\end{Defith}

 The assumptions on $\fctpsi$ imply $\fctpsi(1) \in \alphb$ 
and $\pal_2 = \fctpsi(1) = \delta_1$.
For $\fctpsi(2)$ there are two 
possibilities: either $\fctpsi(2) \in \alphb$ (then $\pal_3 = \fctpsi(1) \fctpsi(2) \fctpsi(1)$ and $\delta_2 = \fctpsi(2)$), or
$\fctpsi(2) = 1$ (then $\pal_3 = \fctpsi(1) \fctpsi(1)$ and
$\delta_2 = \fctpsi(1)$).

Already from this example 
we can see that several functions $\fctpsi$ may lead to the 
same word of first letters $\Delta$: for instance, taking $\fctpsi(2) =
 \fctpsi(1) \in \alphb$ yields the same value of $\delta_2$ as 
taking $\fctpsi(2) = 1 \in \Netoile$, but not the same value of 
$\pal_3$. Using this example it is not difficult to produce 
functions $\fctpsi$ and $\fctpsipr$ with the same
word of first letters but such that $\motdepsi \neq \motdepsipr$.
Therefore a word $\motdepsi$ with \app\ is not given just by
its word of first letters $\Delta$, but by a richer structure: 
the function\footnote{Actually one may restrict to reduced functions,
see \S \ref{subsecreduced} below.} $\fctpsi$.
To be precise, $\fctpsi$ is given exactly by the word $\Delta = \delta_1 \delta_2 \ldots$ together
with the choice, for any $n \geq 1$, of an integer $n' \in \{0,\ldots,n-1\}$ that satisfies either 
$n'=0$ or $\delta_{n'} = \delta_n$. 
If we fix $\Delta$, then a special choice of $\fctpsi$ is obtained by taking for $n'$ the greatest integer
$n' < n$ such that $\delta_{n'} = \delta_n$ (and $n' =0$ if there is no such integer, i.e., if the
letter $\delta_n$ occurs for the first time in $\Delta$ at the $n$-th position). For this function 
$\fctpsi$, the word $\motdepsi$ is the standard
episturmian word with directive word $\Delta$ (see \S \ref{subsecconstrepistu}). 
Therefore 
Definitions \ref{defigalepisu} and \ref{defidirecfctpsi}
generalize Definition \ref{defiepistu} of standard episturmian words.

\bigskip 

\begin{Remark} \label{remunicitetriv}
Two distinct functions $\fctpsi$ and $\fctpsipr$ 
always lead to distinct sequences $(\pal_i)$ and $(\pal'_i)$, but 
may lead to the same word $\motdepsi = \motdepsipr$ (see 
Example \ref{exnecreduit} below).
\end{Remark}

\begin{Example} \label{experio}
If $\fctpsi(n) = n-1$ for any $n \geq N$ then $\pal_{N+\ell} = \pal_N \omega^\ell$ for any $\ell \geq 0$, 
with $\omega = \pal_{N-1}^{-1} \pal_N$. Therefore in this case
$\motdepsi$ is ultimately 
periodic, hence periodic with a palindromic period (see Lemma \ref{lemmewordperiodique} below).
\end{Example}

\begin{Example} \label{exstuepistugal} 
Let $\alphb = \{a,b\}$ be a two-letter alphabet, and $(s_k)_{k \geq 1}$
be a sequence of positive integers.
For any $k\geq 1$, let 
$t_k = s_1+\ldots+s_k$ if $s_1 \geq 2$ and
$t_k = s_1+\ldots+s_{k+1}$ if $s_1 =1$. In both cases, let
$t_0 = 1$. Moreover, let $\fctpsi(i) = i-1$
if $i \geq 1$ is not among $t_0,t_1,t_2,\ldots$, and
 $\fctpsi(t_k) = t_{k-1}-1$ for any $k \geq 2$.
If $s_1 \geq 2$, let
$\fctpsi(1) = a$ and $\fctpsi(t_1) = b$; 
if $s_1 = 1$, let
$\fctpsi(1) = b$ and $\fctpsi(t_1) = a$.
Then the word $\motdepsi$ associated with $\fctpsi$
is the characteristic
Sturmian word with slope $[0,s_1,s_2,\ldots]$. The function $\fctpsi$,
the palindromes $\pal_i$ and the sequence $(t_k)$ are exactly
the same as in \S \ref{subsecstu} (except that the index $k$ in
$t_k$ is shifted if $s_1=1$). 
\end{Example}

\begin{Example} \label{exfiboepistugal}
In the previous example, if $s_k=1$ for any $k \geq 1$ then 
$\fctpsi(1) = b$, $\fctpsi(2) = a$ and $\fctpsi(i) = i-2$ for any $i \geq 3$.
The word $\motdepsi = babbab \ldots$ is the Fibonacci word.
\end{Example}

\subsection{Reduced Functions} \label{subsecreduced}

Two problems immediately  arise
from the construction of words with \app . First, are there 
other palindromic prefixes of $\motdepsi$ than the $\pal_i$'s~? Second, can 
 two distinct functions $\fctpsi$ and $\fctpsi'$ lead to the same word $w$~?

In general, the answers to both questions are positive, as shown in the 
following example. This is the reason why {\em reduced} functions are studied below.

\begin{Example} \label{exnecreduit}
Let $\fctpsi$ be a function, and $i \geq 2$ be a integer, such that
$\fctpsi(i+1) = \fctpsi(i) = i-1$. Let $b_i$ be the finite non-empty word such that
$\pal_i = \pal_{i-1} b_i$. Then $\pal_{i+1} = \pal_{i-1} b_i^2$ and
 $\pal_{i+2} = \pal_{i-1} b_i^4$. Now Lemma \ref{lemfacile} stated below shows that
 $ \pal_{i-1} b_i^3$ is a palindromic prefix of $\pal_{i+2}$ (hence of $\motdepsi$), of length 
 strictly between those of $\pal_{i+1}$ and $\pal_{i+2}$. This gives a palindromic prefix of
 $\motdepsi$ which is not among the $\pal_n$'s constructed from $\fctpsi$. To avoid this problem,
 consider a function $\fctpsipr$ such that $\fctpsipr(n) = \fctpsi(n)$ for $n \leq i$, $\fctpsipr(i+1) = i$
 and $\fctpsipr(i+2) = i+1$. Denoting by $(\pal'_n)$ the sequence of finite palindromes associated with
 $\fctpsipr$, we have $\pal'_n = \pal_n$ for $n \leq i+1$, $\pal'_{i+2} = \pal_{i-1} b_i ^3$
 and $\pal'_{i+3} = \pal_{i-1} b_i ^4$. For $n \geq i+3$, we let $\fctpsipr(n) = \fctpsi(n-1)$ if
 $\fctpsi(n-1) \leq i+1$, and $\fctpsipr(n) = \fctpsi(n-1)+1$ otherwise. Then 
 we have $\pal'_n = \pal_{n-1}$ for any $n \geq i+3$, and $\motdepsi = \motdepsipr$. In this way 
 the functions $\fctpsi$ and $\fctpsipr$ define the same word, but the family of finite palindromes
 associated with $\fctpsipr$ contains the ``missing'' palindrome 
 $\pal_{i-1} b_i^3$.
 \end{Example}

\bigskip

Let 
$\fctpsi: \Netoile \rightarrow \Netoile \undisj \alphb$ be any function (in the sequel
we always assume that, for each $n \geq 1$, either 
$\fctpsi (n) \in \alphb$ or $1 \leq \fctpsi(n) \leq n-1$).

Denote by $(t_k)_{k \geq 0}$ the family of all indexes $n$ (in increasing order) such that 
either $1 \leq \fctpsi(n) \leq n-2$ or $\fctpsi(n) \in \alphb$. This family can be either finite or infinite. 
We always have $t_0 = 1$, since $\fctpsi(1) \in \alphb$.

\begin{Defith} \label{defireduced}
A function $\fctpsi$ is said to be {\em reduced} if the associated sequence 
$(t_k)$ satisfies, for
any $k \geq 1$, the following two conditions:
\begin{itemize}
\item $\fctpsi(t_k) \neq \fctpsi(t_{k-1})$.
\item Either $\fctpsi(t_k) \in \alphb$ or $\fctpsi(t_k) < t_{k-1}$.
\end{itemize}
\end{Defith}

In the special case where the family $(t_k)$ is finite (i.e., 
$\fctpsi(n) = n-1$ for $n$ sufficiently large, see Example \ref{experio}), we assume in this definition
that both properties hold for any $k$ such that $t_k$ exists.

\begin{Remark} The function $\fctpsi$ in Example \ref{exstuepistugal} 
is reduced, and the definition of $(t_k)$ given there is
consistent with the one introduced here.
\end{Remark}

\begin{Remark}
The function $\fctpsi$ in Example \ref{exnecreduit}
is not reduced. In fact there is an integer $k$ such that $i+1 = t_k$, and we have
$t_{k-1} \leq i-1 = \fctpsi(t_k)$. 
\end{Remark}

In the situation of Example \ref{exnecreduit}, we have seen that 
$\fctpsi$ is not reduced, and that the $\pal_i$'s are not the only
palindromic prefixes of $\motdepsi$. Actually both phenomena are equivalent:

\begin{Th} \label{aucunpalmanque} Let $\fctpsi: \Netoile \rightarrow \Netoile \undisj \alphb$ be 
a function such that, for each $n \geq 1$,
either $\fctpsi (n) \in \alphb$ or $1 \leq \fctpsi(n) \leq n-1$. Then the following assertions are equivalent:
\begin{itemize}
\item The function $\fctpsi$ is reduced.
\item The palindromic prefixes of $\motdepsi$ are exactly the $\pal_i$'s
constructed from $\fctpsi$.
\end{itemize}
\end{Th}

This theorem will be proved in the next Section (\S \ref{subsecpreuveaucunpalmanque}). 
It is not difficult
to deduce the following Corollary (see Example \ref{experio} and 
Lemma \ref{lemmewordperiodique}).

\begin{Cor} \label{corapresaucunpm}
Let $\fctpsi$ be a reduced function. Then $\motdepsi$ is
periodic if, and only if, $\fctpsi(n) = n-1$ for any
sufficiently large integer $n$.
\end{Cor}

Let $\fctpsi$ be a reduced function, and $(\pal_i)$ be the associated 
sequence of finite palindromes (that is, thanks to 
Theorem \ref{aucunpalmanque}, the sequence of all palindromic prefixes of
$\motdepsi$). Then the following assertions are easily seen to be
equivalent:
\begin{itemize}
\item For any sufficiently large $i$ we have $\fctpsi(i) \in \Netoile$.
\item For any sufficiently large $i$ we have $\lgr{\pal_{i+1}} \leq
2 \lgr{\pal_i} $.
\end{itemize}
If these assertions hold then $\motdepsi$ can be written on
a finite alphabet.

\bigskip

\bigskip

In addition to Theorem \ref{aucunpalmanque}, another important property of
reduced functions is the following generalization of Theorem 
\ref{thcarac}, proved in Section \ref{secpreuve} below.

\begin{Th} \label{interetreduit}
Let $w$ be an infinite word, and $(n_i)_{i \geq 1}$ 
be the increasing sequence (assumed to be infinite) of the lengths of
its palindromic prefixes (with $n_1 = 0$). Then the following statements are 
equivalent:
\begin{enumerate}
\item[$(i)$] We have $n_{i+1} \leq 2 n_i + 1$ for any $i \geq 1$ (i.e., $w$ has \app ).
\item[$(ii)$] There exists a function $\fctpsi$ such that $w = \motdepsi$.
\item[$(iii)$] There exists a reduced function $\fctpsi$ such that $w = \motdepsi$.
\end{enumerate}
Moreover the reduced function $\fctpsi$ in $(iii)$ is unique.
\end{Th}

It is possible to write down a
``reduction'' algorithm
(generalizing Example \ref{exnecreduit}) that allows one to
obtain, from any function $\fctpsi$, the reduced function $\fctpsipr$ 
such that $\motdepsi = \motdepsipr$. In this situation, the 
 construction of \S \ref{subsecconstruapp} applied with
$\fctpsi$ gives a sequence $(\pal_i)$ of palindromic prefixes of $\motdepsi$; with $\fctpsipr$, it gives another 
sequence $(\pal'_i)$. Theorem \ref{aucunpalmanque} shows that $(\pal_i)$ is a sub-sequence of $(\pal'_i)$.
Again, the ``reduction'' algorithm allows one to obtain explicitly the full sequence
$(\pal'_i)$ from the sub-sequence $(\pal_i)$. 
This algorithm is partly used in \cite{SFOttdio}, but in the present text
we shall not need it; the crucial
point here is just the uniqueness of the reduced function $\fctpsipr$ corresponding to $\fctpsi$.

\begin{Defith} \label{defidirectivefct}
Let $w$ be a word with \app . The reduced function $\fctpsi$ in 
Theorem \ref{interetreduit} is called the {\em directive function} of $w$.
\end{Defith}

\begin{Remark} \label{remunicite}
The uniqueness assertion 
in Theorem \ref{interetreduit} immediately follows from Theorem \ref{aucunpalmanque}
and Remark \ref{remunicitetriv}.
\end{Remark}

Now we can put Definitions \ref{defidirecfctpsi} and \ref{defidirectivefct} 
together in the following way:

\begin{Defith} \label{defidirecmot}
Let $w$ be a word with \app . We call 
{\em word of first letters} associated with $w$ the 
word of first letters associated with the directive function of $w$.
\end{Defith}

The following property holds: if $(\pal_i)$ is the sequence of all palindromic prefixes of a word $w$ with \app , and $\Delta = \delta_1 \delta_2 \ldots$ is the
associated word of first letters, then $\pal_i \delta_i$ is a prefix of
$\pal_{i+1}$ for any $i \geq 1$.

\subsection{Computation of $\densi(w)$ using Reduced Functions} \label{subsecdensipsi}

\begin{Defith} \label{definidepsi}
With any reduced function $\fctpsi$ we associate 
the increasing sequence of non-negative integers $(n_i)_{i \geq 1}$ defined by
 $n_1 = 0$ and, for all $i \geq 1$:
$$n_{i+1} = 2 n_i - n_{\fctpsi(i)} \mbox{ if } \fctpsi(i) \in \Netoile$$
and
$$n_{i+1} = 2 n_i + 1 \mbox{ if } \fctpsi(i) \in \alphb. $$
\end{Defith}

Theorem \ref{aucunpalmanque} shows that 
 $n_i$ is the length of the $i$-th palindromic prefix of $\motdepsi$.
In the same way, we introduce the following definition so that
$\densi(\fctpsi) = \densi(\motdepsi)$:

\begin{Defith} \label{defidensipsi}
For any reduced function $\fctpsi$ we let 
$\densi(\fctpsi) = \limsup \frac{n_{i+1}}{n_i} $, where $(n_i)$ is associated
 with $\fctpsi$
as in Definition \ref{definidepsi}.
\end{Defith}

This definition of $\densi(\fctpsi)$ is completely elementary. It is useful because 
of the following fact: for a word $w$ with \app , we have $\densi(w) = \densi(\fctpsi)$ where 
$\fctpsi$ is the directive function of $w$ (see Definition \ref{defidirectivefct}).

\section{Proof of the Main Results} \label{secpreuve}

\subsection{General Lemmas about Palindromic Prefixes} \label{subseclemmesgeneraux}

The first lemma is very easy, and sufficient
to prove half of Theorem \ref{aucunpalmanque} (see 
\S \ref{subsecpreuveaucunpalmanque} below).

\begin{Lemme} \label{lemfacile} Let $p$ and $u$ be two words, such that
 $p$ and $pu$ are palindromes. Then
$pu^2$ is a palindrome (and so is, by induction, the word
 $pu^n$ for any $n \geq 2$).
Similarly, if $p$ and $up$ are palindromes then $u^np$ is a palindrome for any $n \geq 0$.
\end{Lemme}

\medskip

\Demdeuxpoints If $p$ and $pu$ are palindromes
then we have $\miroir{p} = p$ and $\miroir{u}p = pu$ hence
$$\miroir{pu^2} = \miroir{u} \miroir{u} p = \miroir{u} p u = p u^2$$
The case where $p$ and $up$ are palindromes is analogous.
This concludes the proof of Lemma \ref{lemfacile}.\qed 

\bigskip

In particular, in this situation $pu$ and $pu^2$ are palindromes, one is a prefix of the other, and
the quotient of their lengths is less than 2 (or equal to 2 when $p$ is empty). 
The following lemma gives a kind of converse to this phenomenon (at least in the case $n' = n$).

\begin{Lemme} \label{lemlindepcro}
Let $w$ be an infinite word, and
$n$, $n'$, $n''$ be integers such that
$n' \leq n'' \leq n + n'$. We assume that the prefixes of $w$ with lengths
$n$, $n'$, $n''$ are palindromes, denoted by 
$a$, $a'$ and $a''$ respectively. Let $a_0$ be the prefix of $w$ of length $n+n'-n''$. 
Then the following holds:
\begin{itemize} 
\item There is a word $b$ such that $a = a_0 b$ and $a'' = a' b$.
\item If $n'' \geq n - n'$ then $a_0$ is a palindrome.
\end{itemize}
\end{Lemme}

\begin{Remark} This lemma will be used only when $n'' \geq n - n'$, and in this case the first property will be written
$$a'' = a' a_0 ^{-1} a$$
since $a_0$ is both a suffix of $a'$ and a prefix of $a$.
Moreover, an important special case is when $n = n'$. The lemma then reads: if $a$ and $a'$ are palindromes,
with $n \leq n'' \leq 2n$, then $a_0$ is a palindrome and we have $a = a_0 b$ and $a'' = a_0 b^2$.
\end{Remark}

\medskip

\Dem of Lemma \ref{lemlindepcro}: As $a'$ is a prefix of $a''$, there exists a word $b$ such that
$a'' = a' b$. The word $b$ is a suffix of $a''$, therefore its mirror image
$\miroir{b}$ is a prefix of
$a''$ (hence also of $w$) since $a''$ is a palindrome. Now $\miroir{b}$ 
has length $n'' - n' \leq n$, therefore $\miroir{b}$ 
is a prefix of $a$. As $a$ is a
palindrome, $b$ is a suffix of $a$: there exists a word $c$ such that 
$a = cb$. It is clear that $c = a_0$ is the 
prefix of $w$ of length $ n+n'-n''$. 

Assume now $n'' \geq n - n'$, and let us show that 
$a_0$ is a palindrome. Let $1 \leq i \leq
(n+n'-n'')/2$ ; then we have $i \leq n'$ hence:
$$w_{n+n'-n''+1-i} = w_{n''-n'+i} = w_{n'+1-i} = w_{i}, $$
by using successively that $a$, $a''$ and $a'$
are palindromes. This concludes the proof of Lemma \ref{lemlindepcro}.\qed

\bigskip

\begin{Lemme} \label{lemtrois}
Let $w$ be an infinite word. Let $n' < n''$ be two consecutive lengths of palindromic prefixes of $w$; let us denote by $\pal'$ and $\pal''$ the corresponding prefixes, with
 $\pal'' = \pal' \omega$ for some word $\omega$.
Then any palindromic prefix $\pal$ of $w$ such that $n' \leq \lgr{\pal} \leq n'+n''$
can be written $\pal' \omega^t$ with $t \geq 0$.
\end{Lemme}

\Demdeuxpoints Assume there is a prefix $\pal$ of $w$, of length $n$, which contradicts the lemma and has
minimal length. As $n'$ and $n''$ are consecutive, we have $n > n''$. Lemma \ref{lemlindepcro} gives a 
palindromic prefix $\pal_0$ of $w$ of length $n-\lgr{\omega} > n'$, 
such that $\pal = \pal_0 \omega$. This contradicts the minimality of $\pal$, and 
concludes the proof.\qed

\bigskip

\begin{Lemme} \label{lemsomme} 
Let $w$ be an infinite word. Let $n_0 < n_1 < n_2$ be three consecutive lengths of palindromic prefixes of $w$; let us denote by $\pal_0$, $\pal_1$ and $\pal_2$ the corresponding prefixes. Then:
\begin{itemize}
\item Either $\pal_2 = \pal_1 \pal_0^{-1} \pal_1$,
\item Or $n_2 > n_0 + n_1$.
\end{itemize}
\end{Lemme}

\Demdeuxpoints If $n_2 \leq n_0 + n_1$, one may apply Lemma \ref{lemlindepcro} with 
$n = n_2$, $n' = n_0$ and $n'' = n_1$. Then $n_2 + n_0 - n_1$ is the length of a palindromic prefix of 
$w$; but this length is strictly between $n_0$ and $n_2$, therefore it is $n_1$. We get in this way
$\pal_2 = \pal_1 \pal_0^{-1} \pal_1$, which concludes the proof of the lemma.\qed

\subsection{Ultimately Periodic Words} \label{subsecperiodq} 

\begin{Lemme} \label{lemmewordperiodique}
Let $w$ be an infinite ultimately periodic word, infinitely many prefixes of which are
palindromic. Then $w$ is periodic with a palindromic period. Moreover, 
if $\lgperiw$ denotes the smallest length of a period of
$w$ then there exists $r \in \unlgperiw$ with the following property. For any $n \geq \lgperiw$, the prefix 
of $w$ of length $n$ is a palindrome if, and only if, $n \congru r \mod \lgperiw$.
\end{Lemme}

\Demdeuxpoints If $w$ were ultimately periodic but not periodic, there would exist two
non-empty words 
$\periw_0$ and $\periw$ such that $w = \periw_0 \periw\periw\periw \ldots$, and such that the last letter
of $\periw_0$ be different from that of $\periw$. But this contradicts the assumption that
$w$ has arbitrary long palindromic prefixes. In fact, if we denote by $z_1\ldots z_\lgperiw$ the word $\periw$ and 
by $z_0 \neq z_\lgperiw$ the last letter of $\periw_0$, then this assumption implies that the word
$z_{\lgperiw -1} \ldots z_0$ appears infinitely many times in $w$, and is therefore
a cyclic permutation of the period $z_1\ldots z_\lgperiw$. As $z_0 \neq z_{\lgperiw}$, this is impossible.

Therefore $w$ is periodic, and can be written $w = \periw\periw\periw \ldots $
with a period $\periw$ of minimal length $\lgperiw$. 
Let $n \geq \lgperiw$ be the length of a palindromic prefix of $w$.
Then we have $w_i = w_{n+1-i}$ for all $i \in \unlgperiw$.
If $n'$ is another such integer, not congruent to $n$ mod $\lgperiw$, we obtain
$w_i = w_{i + \eps}$ for all $i \in \unlgperiw$ with $1 \leq \eps \leq \lgperiw - 1$; this contradicts
the minimality of $\lgperiw$. Therefore all lengths of palindromic prefixes lie in the same
congruence class mod $d$; conversely it is clear that any 
$n \geq d$ that belongs to this class is the length of a palindromic prefix of $w$.\qed

\medskip

\begin{Example} For the word $a(abba)^\ome =(aabb)^\ome$, we have $\lgperiw = 4$ and $r = 2$.
\end{Example}

\subsection{Proof of Theorem \ref{aucunpalmanque}} \label{subsecpreuveaucunpalmanque}

Throughout the proof, we fix a function $\fctpsi$, and consider the palindromes $\pal_i$
used to define $\motdepsi$. For any $i \geq 1$
we let $n_i = \lgr{\pal_i}$.

\bigskip

First, let us prove the easier implication (using only Lemma \ref{lemfacile}). Assume 
$\fctpsi$ is not reduced, and all palindromic prefixes of $\motdepsi$ are among the $\pal_i$'s.
There is an index $k \geq 1$ such that either $\fctpsi(t_k) \geq t_{k-1}$ (with 
$\fctpsi(t_k) \in \Netoile$) or $\fctpsi(t_k) = \fctpsi(t_{k-1})$.

Let $j = t_{k-1}$ and $\jpri = t_k$. There are non-empty words $b$ and $b'$ such that
$\pal_{j+1} = \pal_j b$ and 
$\pal_{\jpri+1} = \pal_{\jpri} b'$. Since $\fctpsi(j+1) = \ldots = \fctpsi(\jpri-1) = 1$, we have
$\pal_{\ell} = \pal_j b^{\ell-j}$ for any $\ell \in \{j, \ldots, \jpri\}$ and in particular
$\pal_{\jpri} = \pal_j b^{\jpri-j}$. Applying Lemma \ref{lemfacile} to the palindromes
$ \pal_j b^{\jpri-j-1}$ and $\pal_{\jpri}$ proves that $\pal_{\jpri} b$ is a palindrome. 
If $b' = b^n$ with $n \geq 2$, this is a palindromic prefix of $\pal_{\jpri+1}$ (hence of $\motdepsi$), whose length
is strictly between those of $\pal_{\jpri}$ and $\pal_{\jpri+1}$; this contradicts the assumption that all
 palindromic prefixes of $\motdepsi$ are among the $\pal_i$'s.
Therefore $b'$ cannot be a non-trivial power of $b$.

In the case where $\fctpsi(\jpri) \in \Netoile$ and $\fctpsi(\jpri) \geq j$, we have
$b' = b^{\jpri-\fctpsi(\jpri)}$ with $\fctpsi(\jpri) \leq \jpri-2$, hence a contradiction.

Assume now $\fctpsi(j) = \fctpsi(\jpri) \in \alphb$. Then $b = \fctpsi(j) \pal_j$ hence
$b' = \fctpsi(j) \pal_{\jpri} = \fctpsi(j) \pal_j b^{\jpri-j} = b^{\jpri-j+1}$ with 
$\jpri-j+1 \geq 2$; this is again a contradiction.

At last, assume $\fctpsi(j) = \fctpsi(\jpri) \in \Netoile$. Then 
we have in the same way
$b = \pal_{\fctpsi(j)}^{-1} \pal_j$ hence
$b' = \pal_{\fctpsi(j)}^{-1} \pal_{\jpri} = b^{\jpri-j+1}$, hence a contradiction.
This concludes the proof of the first implication in Theorem \ref{aucunpalmanque}.

\medskip

Let us prove the converse now.
Assume $\fctpsi$ is reduced, and let $w = \motdepsi$. 
Let $\pal''$ be the palindromic prefix of $w$ of minimal length among those which are not $\pal_i$'s.
Let $i$ be the integer such that 
 $\lgr{\pal_i} < \lgr{\pal''} < \lgr{\pal_{i+1}}$.
 Since $\pal_1 = \motvide$ and $\pal_2 = \fctpsi(1) \in \alphb$, we have $i \geq 2$.

Let $\omega$ be such that $\pal'' = \pal_i \omega$. 
Then Lemma \ref{lemtrois} gives an integer $t \geq 2$ 
such that $\pal_{i+1} = \pal_i \omega^t$; the definition of $\pal_{i+1}$ shows that
$\omega$ is a suffix of $\pal_i$. Now Lemma \ref{lemlindepcro} (with $n=n'=n_i$ and $n'' = \lgr{\pal''}$)
implies that
$ \pal_i \omega^{-1}$ is a palindromic prefix of $w$.

Let us prove that $i$ is among the $t_k$'s. This is obvious if $\fctpsi(i) \in \alphb$, so we may assume
 $\fctpsi(i) \in \Netoile$. Then $\pal_{i} = \pal_{\fctpsi(i)} \omega^t$ and the palindromic prefix
$ \pal_i \omega^{-1}$ has length strictly between
$\lgr{\pal_{\fctpsi(i)}}$ and $\lgr{\pal_i}$ since $t \geq 2$. By minimality of $\pal''$, this implies
$\fctpsi(i) \leq i-2$, and concludes the proof that 
there exists $k \geq 1$ such that $i = t_k$. 

Let $j = t_{k-1}$, and $b$ be the word such that $\pal_{j+1} = \pal_{j} b$. Then we have $\pal_i = \pal_j b^{i-j}$,
since $\fctpsi(\ell) = \ell - 1$ for any $\ell \in \{j+1,\ldots, i-1\}$. 
Now we have $t \geq 2$, hence $2 \lgr{\omega} \leq n_i + 1$ and Lemma \ref{lemsomme} yields:
$$n_i - \lgr{\omega} \geq \frac{n_i-1}{2} \geq \frac{n_j+\lgr{b}-1}{2} \geq \frac{n_j+n_{j-1}}{2} > n_{j-1}.$$
Therefore the palindromic prefix $ \pal_i \omega^{-1}$ is among $\pal_j$, \ldots, $\pal_{i-1}$. This shows that
$\omega$ is a power of $b$, say $\omega = b^u$.

\bigskip

First, let us assume $\fctpsi(i) \in \alphb$. If we also have
$\fctpsi(j) \in \alphb$ then $b = \fctpsi(j)\pal_j$ and
$\fctpsi(i)\pal_i = \omega^t = b^{ut} = (\fctpsi(j)\pal_j)^{ut}$ hence
$\fctpsi(i) =\fctpsi(j) \in \alphb$, which is a contradiction. Now in the case 
$\fctpsi(j) \not\in \alphb$ we have
$0 \leq n_{\fctpsi(j)} < n_{j-1} \leq \lgr{b} - 1$ by Lemma \ref{lemsomme}. However
$n_{\fctpsi(j)} \congru n_j \congru n_i \congru -1 \mod \lgr{b}$ since 
$\fctpsi(i) \pal_i = \omega^t$; 
this is again a contradiction. 

\bigskip

To conclude the proof, we have to consider the case where $\fctpsi(i) \in \Netoile$.
Since $\fctpsi$ is reduced, this gives
$\fctpsi(i) < t_{k-1}=j$. We see that
$n_{\fctpsi(i)} = n_i - t \lgr{\omega}$ belongs to $n_i + \lgr{b} \Z = n_j + \lgr{b} \Z$.
Now $n_j \leq 2 n_{j-1} + 1 < 2 \lgr{b}$ thanks to Lemma \ref{lemsomme}, hence 
$n_{\fctpsi(i)} = n_j - \lgr{b}$. 
This implies $\lgr{b} \leq n_j$, hence $\fctpsi(j) \not\in \alphb$. But we then have
$ n_{\fctpsi(j)} = n_j - \lgr{b} = n_{\fctpsi(i)}$, and this equality contradicts the
assumption that $\fctpsi$ is reduced. 

\bigskip

This concludes the proof of Theorem \ref{aucunpalmanque}.\qed

\subsection{Proof of Theorem \ref{interetreduit}} \label{subsecdemdeuxth}
 
The uniqueness statement at the end of Theorem \ref{interetreduit} follows from
Theorem \ref{aucunpalmanque} and Remark \ref{remunicitetriv} (as noticed in Remark \ref{remunicite}
above). Let us prove that $(i)$, $(ii)$ and $(iii)$ in Theorem \ref{interetreduit}
are equivalent.

The implications $(iii) \doublerightarrow 
(ii)$ and 
 $(ii) \doublerightarrow (i)$ are obvious; let us prove that $(i)$ implies $(iii)$. Assume that 
$w$ satisfies $n_{i+1} \leq 2 n_i + 1$ for all $i \geq 1$. We denote by $\pal_i$ the palindromic prefix of $w$ with length $n_i$.

For any $i \geq 1$, let us define $\fctpsi(i)$ in the following way. If $n_{i+1} =2 n_i +1$, we let
$\fctpsi(i) \in \alphb$ be the letter in $\pal_{i+1}$ that comes right after $\pal_i$ (that is, the central letter
of the palindrome $\pal_{i+1}$ which has odd length). Otherwise, we apply Lemma \ref{lemlindepcro} with
$n=n'=n_i$ and $n'' = n_{i+1}$. This gives an integer $\fctpsi(i)$ between 1 and $i-1$ such that 
$\pal_{i+1} = \pal_i \pal_{\fctpsi(i)}^{-1} \pal_i$. 

With this construction, it is clear that $w = \motdepsi$. More precisely, 
the palindromes $\pal_i$ are
exactly those constructed using $\fctpsi$ in \S \ref{subsecconstruapp}. 
Since they are (by hypothesis)
the only palindromic prefixes of $w$, Theorem \ref{aucunpalmanque} proves that $\fctpsi$ is reduced.

This concludes the proof of Theorem \ref{interetreduit}.\qed

\section{Palindromic Prefix Density} \label{secdensite}

In \S \ref{subsecasyabondant}, we show how the results on words with \app\ 
can be generalized to words such that
$n_{i+1} \leq 2n_i + 1$ for any sufficiently large $i$. This enables us in 
\S \ref{subsecdemtheoqcq} to describe $\specqcq$ in terms of $\densi(\fctpsi)$ 
for reduced functions $\fctpsi$, and to prove Theorem
 \ref{theoqcq} stated in the Introduction. This description
 makes use of a technical statement (Proposition \ref{propindep}), the proof of which
 is postponed to the 
Appendix (Section \ref{secappendice}).

\smallskip

Given an infinite word $w$, we denote by 
$(n_i)$ the increasing sequence of all lengths of
 palindromic prefixes of $w$ (with $n_1 = 0$). 
We let $\pal_i$ be the palindromic prefix of $w$ with length $n_i$.

\subsection{Words with Asymptotically Abundant Palindromic Prefixes} \label{subsecasyabondant}

The following proposition is a generalization of results stated in Section \ref{secapp}. It
enables one to construct all words satisfying 
$n_{i+1} \leq 2n_i + 1$ for any sufficiently large $i$.

\begin{Prop} \label{propoutilqcq}
Let $w$ be a word with infinitely many palindromic prefixes, 
and $\imin$ be an integer. The following statements are equivalent:
\begin{enumerate}
\item[$(i)$] We have $n_{i+1} \leq 2n_i + 1$ for any $i \geq \imin$.
\item[$(ii)$] There exists a reduced function $\fctpsi$ such that, for any $i \geq \imin$: 
\begin{itemize}
\item Either $\fctpsi(i) \in \alphb$ and $\pal_{i+1} = \pal_i \fctpsi(i)\pal_i $,
\item Or $\fctpsi(i) \in \Netoile$ and $\pal_{i+1} = \pal_i \pal_{\fctpsi(i)}^{-1} \pal_i$.
\end{itemize}
\end{enumerate}
\end{Prop}

This proposition means that the palindromic prefixes $\pal_i$
of $w$ satisfy (for $i$ sufficiently large) the same recurrence relation
as those of the word with \app\ $\motdepsi$. This recurrence
relation is completely determined by the function $\fctpsi$.

Let $\fctpsi$ be a reduced function, and $\pal_{\imin}$ be a finite word with exactly
$\imin$ palindromic prefixes (including $\pal_1 = \motvide$ and $\pal_{\imin}$ itself),
denoted by $\pal_1$, \ldots, $\pal_{\imin}$.
Then using the construction of \S \ref{subsecconstruapp} 
(which is the same as $(ii)$ in Proposition \ref{propoutilqcq}) for $i \geq \imin$ 
one obtains an infinite word
$w$ with infinitely many palindromic prefixes, 
of which $\pal_{\imin}$ is a palindromic prefix. For this
word we have
$n_{i+1} \leq 2n_i + 1$ for any $i \geq \imin$, but there is no reason why this relation
would hold for $i < \imin$: for instance if $\imin = 2$,  $\pal_{\imin}$ 
may  be any  palindrome whose only palindromic prefixes are $\motvide$, its first letter and itself, so it might have length greater than 3.  
Moreover, Proposition \ref{propoutilqcq} implies that any word $w$ satisfying $(i)$ can be obtained 
in this way. 

\bigskip

\Dem of Proposition \ref{propoutilqcq}: 
The implication $(ii) \doublerightarrow (i)$ is clear. To prove the converse, 
we define $\fctpsi(i)$ for $i \geq \imin$ as in the proof of 
Theorem \ref{interetreduit} in \S \ref{subsecdemdeuxth}. For $i < \imin$, we 
let $\fctpsi(i) = a$ if $i$ is even, and $\fctpsi(i) = b$ if $i$ is odd, where $a \neq b$ are 
two elements in $\alphb$ (the trivial case where the alphabet contains only one letter is easily
dealt with). All assertions in Proposition \ref{propoutilqcq} immediately follow, except the fact
that $\fctpsi$ is reduced. To prove this fact, one follows exactly the same lines as in the proof of
Theorem \ref{aucunpalmanque} in \S \ref{subsecpreuveaucunpalmanque}.\qed

\subsection{Elementary Description of $\specqcq$ and $\spectre$} \label{subsecdemtheoqcq}

To establish a relationship between $\specqcq$ and $\spectre$, we shall use
the definitions and statements of \S \ref{subsecdensipsi}, and
 the following (technical) proposition proved in the Appendix (\S \ref{subsecpreuvepropindep}):

\begin{Prop} \label{propindep}
Let $\fctpsi$ be a reduced function such that $\densi(\fctpsi) < 2$, $(n_i)$ be the associated sequence, and 
 $(n'_i)_{i \geq 1}$ be another increasing sequence of non-negative integers such that 
 $n'_{i+1} = 2 n'_i - n'_{\fctpsi(i)}$ for $i$ sufficiently large. 

Then the quotient $n'_i / n_i$ has a finite positive limit as $i$ tends to infinity, and we have:
$$\densi(\fctpsi) = \limsup \frac{n'_{i+1}}{n'_i}.$$ 
\end{Prop}

Let $w$ be any word such that $\densi(w) < 2$, and $(n_i)_{i \geq 0}$ be the sequence
of all lengths of palindromic prefixes of $w$. Proposition \ref{propoutilqcq} yields
a reduced function $\fctpsi$ such that $n_{i+1} = 2n_i - n_{\fctpsi(i)}$ for any sufficiently large
integer $i$. This proves that $(n_i)$ satisfies the same recurrence relation as
the sequence associated with $\fctpsi$
(in Definition \ref{definidepsi}), but the initial values may be distinct.
Proposition \ref{propindep}
shows that these initial values have no influence on $\limsup \frac{n_{i+1}}{n_i}$,
hence $\densi(w) = \densi(\fctpsi)$. This proves that $\specqcq \inter [1,2)$ is contained
in the set of values $\densi(\fctpsi)$ for reduced functions $\fctpsi$, which is
exactly $\spectre$ (see \S \ref{subsecdensipsi}). Since $\spectre$ is obviously contained in
$\specqcq$, this gives $\specqcq \inter [1,2) = \spectre \inter [1,2)$. Now 
considering a characteristic Sturmian word whose slope has unbounded partial 
quotients proves (thanks to Equation \eqref{eqstu}) that
$2 \in \spectre \inclus \specqcq$. Therefore we have proved the following result:

\begin{Th} \label{theoauxintro}
Both $\spectre$ and $\specqcq \inter [1,2]$ are exactly the set of values
taken by $\densi(\fctpsi)$ for reduced functions $\fctpsi$.
\end{Th}

Together with the inclusion $(2 , +\infty] \inclus \specqcq$ proved in 
\S \ref{subsecscarse}, this statement implies Theorem \ref{theoqcq}.

\section{Peculiar Study around  $\sqrt3$} \label{secsqrtt}

In this Section, we focus on the sets $\spectre$ and $\specBL$ 
around $\sqrt3$. First, we prove Theorem \ref{theosqrttrois}
stated in the Introduction, namely 
$\spectre \inter [1, \sqrt3] = \specBL \inter [1, \sqrt3] $ (\S \ref{subsecpreuvetheosqrttrois}).
Next, we define and 
study $\alphb$-strict words with \app\ (\S \ref{subsecinitial}); this allows 
us to prove that any $\alphb$-strict
word $w$ with \app\ such that 
$\densi(w) < \sqrt3$ is either periodic or characteristic Sturmian.

At last, we prove that Theorem \ref{theosqrttrois} is optimal, since
there is no gap in $\spectre$ right above $\sqrt3$ 
(\S \ref{subsecracinetroisnonepistu}) whereas there is one
in $\specBL$ (\S \ref{subsecspectresqrt3}). 

\smallskip

As in the previous Section, we shall use a technical result
(Lemma \ref{lemdecalcul}, stated in \S \ref{subsecpreuvetheosqrttrois}) 
the proof of which is postponed to the Appendix.

\smallskip

Given an infinite word $w$, we still denote by 
$(n_i)$ the increasing sequence of all lengths of
 palindromic prefixes of $w$ (with $n_1 = 0$). 
We let $\pal_i$ be the palindromic prefix of $w$ with length $n_i$.

\subsection{Proof of Theorem \ref{theosqrttrois}} \label{subsecpreuvetheosqrttrois}

Let $\fctpsi$ be a reduced function. As in \S \ref{subsecreduced}, we denote
 by $(t_k)_{k \geq 0}$ the family of all indexes $n$
 (in increasing order) such that 
either $1 \leq \fctpsi(n) \leq n-2$ or $\fctpsi(n) \in \alphb$. 

Assume that $\densi(\fctpsi) < 2$. Then $\fctpsi(i) \in \Netoile$ for any sufficiently large
integer $i$, hence:
\begin{itemize}
\item Either the family $(t_k)$ is finite, that is 
$\fctpsi(i) = i-1$ for any
sufficiently large $i$. This obviously implies $\densi(\fctpsi) = 1$.
\item Or the family $(t_k)$ is infinite, and
for $k$ sufficiently large we have $\fctpsi(t_k) < t_{k-1}$.
\end{itemize}
The first condition corresponds to periodic words (see
Corollary \ref{corapresaucunpm}). In the second
condition, the special case where
$\fctpsi(t_k) = t_{k-1} - 1$ for any $k \geq 2$ corresponds to 
characteristic Sturmian words (see Example \ref{exstuepistugal}). 
The following result shows that any reduced function $\fctpsi$
with $\densi(\fctpsi) < \sqrt3$ is either ``periodic'' or 
``asymptotically Sturmian''.

\begin{Lemme} \label{lemdecalcul}
Let $\fctpsi $ be a reduced function such that $\densi(\fctpsi) < \sqrt3$. Then either the family
$(t_k)$ is finite, or for any $k$ sufficiently large we have
$$\fctpsi(t_k) = t_{k-1} - 1.$$
\end{Lemme}

The value $\sqrt3$ in this lemma is optimal (see the
end of \S \ref{subsecracinetroisnonepistu}). 
We postpone the proof (which is completely elementary, but rather 
technical) to the Appendix (\S \ref{subsecpreuvelemdecalcul}).

\bigskip

Now let us prove 
 Theorem \ref{theosqrttrois} stated in the Introduction, namely 
$\spectre \inter [1, \sqrt3] = \specBL \inter [1, \sqrt3] $. 
 First of all, it is readily seen that 
 $\sqrt3 = [1,1,\overline{2,1}] \in \specBL \inclus \spectre $ (see the beginning
 of \S \ref{subsecspectresqrt3}). Let $\fctpsi$ be a reduced function such that
$\densi(\fctpsi) < \sqrt3$. If $\motdepsi$ is periodic then 
Corollary \ref{corapresaucunpm} implies $\densi(\fctpsi) = 1 \in
\specBL$. Otherwise the associated sequence $(t_k)$ is infinite, 
 and satisfies
$\fctpsi(t_k) = t_{k-1}-1$ for any sufficiently large $k$
thanks to Lemma \ref{lemdecalcul}. Let
$s_k = t_k - t_{k-1}$ for $k \geq 1$, and $\fctpsipr$ be the function defined from
$(s_k)$ in Example \ref{exstuepistugal}. Proposition \ref{propindep} shows 
 that $\densi(\fctpsi) = \densi(\fctpsipr) = \densi(\caral)$,
where $\caral$ is the characteristic Sturmian word with slope
$[0,s_1,s_2,\ldots]$. Thanks to Theorem \ref{theoauxintro}, this 
concludes the proof of Theorem \ref{theosqrttrois}.\qed

\bigskip

Actually the result we have proved is slightly more precise that Theorem \ref{theosqrttrois}:
for any non-periodic word $w$ such that $\densi(w) < \sqrt3$ we have found a
 characteristic Sturmian word $\caral$ such that the palindromic prefixes of $w$ and those
 of $\caral$ satisfy (asymptotically) the same recurrence relation. In the next
 Section, we show that two additional assumptions (namely
abundance of palindromic prefixes and 
 $\alphb$-strictness) imply $w = \caral$.
 
 However there is a characteristic Sturmian word $w$, and 
 a non-episturmian $\alphb$-strict word $w'$ with \app , such that
 $\densi(w) = \densi(w') = \sqrt3$ (see the end of \S \ref{subsecracinetroisnonepistu}).
 This shows that $\densi(\fctpsi)$ does not characterize a reduced function $\fctpsi$.

\subsection{Initial Values and Strict Words} \label{subsecinitial}

Let us consider 
(on the three-letter alphabet $\alphb = \{a,b,c\}$)
the finite word $\pal_4 = abacaba$. It has four palindromic prefixes:
$\pal_1 = \motvide$, $\pal_2 = a$, $\pal_3 = aba$ and $\pal_4$.
It is a palindromic prefix of the word $\motdepsi$ constructed 
(as in \S \ref{subsecconstruapp})
from any function
$\fctpsi$ such that $\fctpsi(1) = a$, $\fctpsi(2) = b$ and $\fctpsi(3) = c$.
Now if $\fctpsi$ is given by $\fctpsi(i) = i-2$ for any $i \geq 4$ then $\motdepsi$ 
behaves ``asymptotically'' like the Fibonacci word (see Example \ref{exfiboepistugal}); for
instance $\densi(\motdepsi) = \gamma$. 

This word $\motdepsi$ is a standard episturmian word, which is not Sturmian (it cannot be written 
on a two-letter alphabet), but which behaves like a Sturmian word. 
To avoid this kind of examples, one usually restricts to $\alphb$-strict 
 standard episturmian words (\cite{DJP}, \S 4.2), also known 
 as characteristic Arnoux-Rauzy words. 
 
 Now let us turn to the (more general) case of words with \app . 
 There are words which have \app , behave 
 like a Sturmian word (as above), but are not Sturmian -- and not episturmian either. We introduce
the following definition (recall that with any
word $w$ with \app\ we associate in Definition \ref{defidirecmot} its
word of first letters $\Delta = \delta_1 \delta_2 \ldots$):

\begin{Defith} A word $w$ on an alphabet $\alphb$, with \app , is said to be $\alphb$-strict
if every letter in $\alphb$ occurs infinitely many times in the word of first letters of $w$.
\end{Defith}

This definition extends that of $\alphb$-strict 
 standard episturmian words (\cite{DJP}, \S 4.2). It allows us to
state the following result (recall from \cite{DJP}, Theorem 4, that 
$\alphb$-strict standard episturmian words on a two-letter alphabet $\alphb$
are exactly characteristic Sturmian words):

\begin{Th} \label{densistrict} Let $w$ be any non-periodic
word with \app . We assume $w$ to be $\alphb$-strict, and 
$\densi(w) < \sqrt3$. Then $\alphb$ contains exactly two letters, and $w$ is characteristic Sturmian.
\end{Th}

As a special case, we get the following result:

\begin{Cor}
 Let $w$ be any non-periodic $\alphb$-strict standard episturmian word, with
$\densi(w) < \sqrt3$. Then $\alphb$ contains exactly two letters, and $w$ is characteristic Sturmian.
\end{Cor}

\Dem of Theorem \ref{densistrict}: Denote by $(n_i)$ the increasing sequence of all lengths
of palindromic prefixes of $w$, and by $\fctpsi$ be the directive function of
$w$. Lemma \ref{lemdecalcul} shows that $\fctpsi(t_k) = t_{k-1}-1$ for
 any sufficiently large $k$. Denote by $\Delta$ the word of first letters associated with
 $\fctpsi$ (i.e., with $w$). Then $\delta_n = \delta_{n-1}$ if $n$ is not among the $t_k$'s,
 and $\delta_{t_k} = \delta_{t_{k-1}-1} = \delta_{t_{k-2}}$ for any 
 sufficiently large $k$. Therefore $\delta_n$ takes infinitely many times at most two
 values. Since $w$ is $\alphb$-strict, $\alphb$ contains at most two letters.
 If $\alphb$ is reduced to a single letter then $w$ is nothing but the periodic repetition
 of this letter. Otherwise $w$ is a standard $\alphb$-strict episturmian word
 on a two-letter alphabet, hence is characteristic Sturmian 
 (\cite{DJP}, Theorem 4).\qed

\subsection{Non-Episturmian Examples near $\sqrt3$} \label{subsecracinetroisnonepistu}

The following proposition, together with Proposition \ref{propstuaudelasqrtt}
proved in \S \ref{subsecspectresqrt3}, shows that the value $\sqrt3$ is optimal in 
Theorem \ref{theosqrttrois}.

\begin{Prop} There exists a decreasing sequence in $\spectre$ with limit $\sqrt3$.
\end{Prop}

To prove this, consider for any integer $n \geq 2$ 
the function $\fctpsi_n$ defined as follows (with the 
two-letter alphabet $\alphb = \{a,b\}$).

Let $\phi_n: \{0, \ldots, 4n\} \rightarrow \{1,2,3\}$ be defined by:
\begin{itemize}
\item $\phi_n(0) = 3$.
\item For $ i \in \{1, \ldots,   n\}$, $\phi_n (i) = 2$.
\item For $ i \in \{n+1, \ldots,   4n\}$ with $i \congru n+1 \mod 3$, $\phi_n(i) = 2$.
\item For  $ i \in \{n+1, \ldots,   4n\}$ with $i \congru n+2 \mod 3$, $\phi_n(i) = 1$.
\item For  $ i \in \{n+1, \ldots,   4n\}$  with $i \congru n \mod 3$, $\phi_n(i) = 3$.
\end{itemize}
Then we let $\fctpsi_n(1) = a$, $\fctpsi_n(2) = b$ and 
$\fctpsi_n(i) = i - \phi_n(i')$ for any $i \geq 3$, where $i'$ is the integer
between 0 and $4n$ which is congruent to $i$ modulo $4n+1$.

The word of first letters associated with $\fctpsi_n$ is (in the case where $n \geq 3$ is odd)
$$\Delta_n = \Big( (ab)^n (bba)^n b \Big) ^\ome .$$
By definition, for any $i \geq 3$ we have $\delta_{\fctpsi_n(i)} = \delta_i$, with
$\Delta_n = \delta_1 \delta_2 \ldots$.
Moreover, if $i \geq 3$ is not a multiple of $4n+1$ then $\fctpsi_n(i)$ is the largest 
index $i'$ such that $\delta_{i'} = \delta_i$. But when $i$ is a multiple of $4n+1$,
we have $\fctpsi_n(i) = i-3$ with $\delta_{i-3} = \delta_{i-2} = 
\delta_i = b$. 
Moreover it is easily checked that $\fctpsi_n$ is reduced. 
Therefore Theorem \ref{aucunpalmanque} shows that $\motdepsin$ is not
a standard episturmian word.

\medskip

It is possible to prove that $\densi(\fctpsi_n)$ is greater than $\sqrt3$, and tends
to $\sqrt3$ as $n$ tends to infinity.\qed

\bigskip

Using the same ideas, it is possible to construct 
a word $ \motdepsi$ (which is not standard episturmian) such that 
$\densi(\motdepsi) = \sqrt3$, by choosing an increasing sequence $(n_k)$, with sufficiently fast growth, and
building a function $\phi$ by concatenating the functions $\phi_{n_k}$.
This proves that the conclusion of Lemma \ref{lemdecalcul} does not hold 
for any functions $\fctpsi$ such that $\densi(\fctpsi) = \sqrt3$.

\begin{Remark} The word  $ \motdepsi$ constructed in this way is not episturmian; however, the value of $\densi(\motdepsi)$, namely $\sqrt3$, equal to $\densi(w')$ for some characteristic Sturmian word $w'$ (see the beginning of \S \ref{subsecspectresqrt3}). This proves that knowing $\densi(w)$ does not provide information on the structure of $w$.
\end{Remark}

\subsection{The Sturmian Spectrum near $\sqrt3$} \label{subsecspectresqrt3}

To prove that
$$\sqrt3 = [1,1,\overline{2,1}] = 1.7320\ldots$$
belongs to $\specBL$, it is enough to apply Equation \eqref{eqstu}
with $s_k = 1$ for $k$ even and $s_k = 2$ for $k$ odd. 
Now taking $s_k = 3$ for any $k$ yields another element of $\specBL$:
$$\frac{7+\sqrt{13}}{6} = [1,1,\overline{3}] = 1.7675\ldots . $$

The following proposition proves that there is nothing inbetween:

\begin{Prop} \label{propstuaudelasqrtt}
There is no element of $\specBL$ between $\sqrt3$ and $\frac{7+\sqrt{13}}{6}$.
\end{Prop}

\Demdeuxpoints For a sequence $b = (b_1, b_2, \ldots)$ of positive integers, and a non-negative integer $k$,
we let $T^k b = (b_{k+1}, b_{k+2}, \ldots)$. Cassaigne proved in \cite{Cassaigne} that $\specBL$ is
the set of numbers $[1,1,b_1,b_2,\ldots]$ where the sequence $b$ satisfies 
$[b] \geq [T^k b]$ for any $k \geq 0$ (where $[b]$ is the continued fraction $[b_1,b_2,\ldots]$). 
Let $b$ be such a sequence, with $\sqrt3 < [1,1,b_1,b_2,\ldots] < \frac{7+\sqrt{13}}{6}$.

First of all, let us prove that $b_i \in \{1,2\}$ for any $i \geq 1$. Indeed, the assumption
$[1,1,b_1,b_2,\ldots] < [1,1,\overline{3}]$ means $[b] < [\overline{3}]$ hence $b_1 \leq 3$. If $b_1 \leq 2$ then
the assertion is proved (since $[b] \geq [T^k b]$ for any $k$). Otherwise $b_1 = 3$ and
$b_i \in \{1,2,3\}$ for any $i$. But $[b_2,b_3,\ldots] > [\overline{3}]$ yields $b_2 \geq 3$ hence $b_2 = 3$. Now
$[b] \geq [T b]$ gives $[3,b_3,\ldots] < [b_3,\ldots,]$ hence $b_3 = 3$. Repeating these arguments gives
$b_i = 3$ for any $i$, in contradiction with the assumption.

Now we have $[b] > [\overline{2,1}]$ with $b_i \in \{1,2\}$ for all $i$. This gives $b_1 = 2$ and
$[b_2,b_3,\ldots] < [\overline{1,2}]$ hence $b_2 = 1$. Repeating this process yields $[b] = [\overline{2,1}]$,
that is a contradiction.

This concludes the proof.\qed

\section{Open Questions} \label{secquestions}

This Section is devoted to open questions about words $w$ with \app . Some of them were asked 
by various specialists I would like to thank.

We let $w$ be a word with \app .

\bigskip

Do letters (and more generally factors) have frequencies in $w$? 
Is it possible to compute these frequencies in terms of the directive function of $w$ (see Definition \ref{defidirectivefct})?

What is the complexity of $w$?

What are the recurrence quotient, and the critical exponent of 
$w$ (in terms of $\fctpsi$)? 

What is the bound one should put instead of $\sqrt3$ in Theorem \ref{densistrict} to ensure that 
$\alphb$ contains at most 3 letters? More generally, for any integer $p$ one could study 
the set of values taken by $\delta(w)$, for words $w$ written on an
alphabet of at most $p$ letters.

Are there other ways to construct the set of all words with \app\ (as 
for standard episturmian words)? 

Is is possible to write words with \app\ as fixed points of morphisms?

Is there a way to define ``non-standard'' words with abundant palindromic prefixes~? 
This would be a class of words that
behaves with respect to words with \app\ in the same way as Sturmian words 
with respect to characteristic Sturmian words, and in the same way as episturmian words with 
respect to standard episturmian words.

What does $\specBL$ look like between the least accumulation point $\sigma_{\infty}$ and
$\sqrt3$ (see \cite{BL}, \S 8)? What does $\spectre$ look like above $\sqrt3$? In which intervals is it dense?

\section{Appendix: Study of Some Recurrence Relations} \label{secappendice}

In this Appendix, we study the linear recurrence relation satisfied by the sequence
$(n_i)$ of all lengths of palindromic prefixes of a word $w$ with \app . We focus on
asymptotic properties of this sequence (and especially on 
$\limsup n_{i+1}/n_i$). The point of view is to forget everything
about words: all statements and proofs are completely elementary, and rely only
on the recurrence relation associated with a reduced function $\fctpsi$. 

The point is to prove two technical statements used in the text: Proposition \ref{propindep}
and Lemma \ref{lemdecalcul}. Moreover, the tools introduced here are used in 
\cite{SFOttdio}.

\subsection{Definitions} \label{subsecappendicedef}

Recall from \S \ref{subsecconstruapp} that we consider functions
$\fctpsi: \Netoile \rightarrow \Netoile \undisj \alphb$ such that, for each $n \geq 1$:
$$\mbox{ either } \fctpsi (n) \in \alphb \mbox{ or } 1 \leq \fctpsi(n) \leq n-1. $$
With such a function we associate in \S \ref{subsecreduced} a (finite or infinite)
family $(t_k)_{k \geq 0}$, namely the family of all indexes $n \geq 1$ (in increasing order) such that 
either $1 \leq \fctpsi(n) \leq n-2$ or $\fctpsi(n) \in \alphb$. 

We recall that $\fctpsi$ is {\em reduced} (see Definition \ref{defireduced}) if, for
any $k \geq 1$, the following two conditions are satisfied:
\begin{itemize}
\item $\fctpsi(t_k) \neq \fctpsi(t_{k-1})$.
\item Either $\fctpsi(t_k) \in \alphb$ or $\fctpsi(t_k) < t_{k-1}$.
\end{itemize}
In the special case where the family $(t_k)$ is finite (i.e.,
$\fctpsi(n) = n-1$ for $n$ sufficiently large), we assume in this definition
that both properties hold for any $k$ such that $t_k$ exists.

Moreover, with any reduced function $\fctpsi$ we associate (as in \S \ref{subsecdensipsi}) 
the increasing sequence of non-negative integers $(n_i)_{i \geq 1}$ defined by
 $n_1 = 0$ and, for all $i \geq 1$:
$$n_{i+1} = 2 n_i - n_{\fctpsi(i)} \mbox{ if } \fctpsi(i) \in \Netoile$$
and
$$n_{i+1} = 2 n_i + 1 \mbox{ if } \fctpsi(i) \in \alphb, $$
and we let 
$$\densi(\fctpsi) = \limsup \frac{n_{i+1}}{n_i}. $$

\subsection{First Lemmas}

\begin{Lemme} \label{lemdelemtri}
Let $\fctpsi$ be a reduced function, and 
$(n_i)_{i \geq 1}$ be the associated sequence.
Then for any $k \geq 1$ we have, with $i = t_k$:
$$ 
n_{i+1} > n_i + n_{i-1} .
$$
\end{Lemme}

\Demdeuxpoints We proceed
by induction on $k$. We have $t_0 = 1$
hence $\fctpsi(t_1) \in \alphb$ since $\fctpsi$ is reduced. Therefore $n_{t_1+1}
=2n_{t_1}+1 > n_{t_1-1}$, which proves the result for $k=1$. 
Assume it holds for $k$. If $\fctpsi(t_{k+1}) \in \alphb$ then it clearly holds for $k+1$; otherwise
we have:
$$n_{t_{k+1}+1} \geq 2 n_{t_{k+1}} - n_{t_k - 1} = n_{t_{k+1}} + n_{t_{k+1} - 1} + n_{t_k+1} - n_{t_k}- n_{t_k - 1}.$$
This concludes the proof.\qed

\bigskip

\begin{Lemme} \label{lembornee}
Let $\fctpsi$ be a reduced function. 
Then $\densi(\fctpsi) < 2$ if, and only if, 
 there is an integer $\borne$ such that, for any sufficiently large $i$:
$$\fctpsi(i) \in \Netoile \mbox{ and } \fctpsi(i) \geq i - \borne. $$
\end{Lemme}

\Demdeuxpoints Denote by $(n_i)$ the sequence associated with $\fctpsi$.

Assume there is an integer 
$\borne$ be such that $\fctpsi(i) \geq i - \borne$ for all $i$ sufficiently large.
Then for $i$ sufficiently large we have 
$\fctpsi(i) \leq 3n_i$ hence
 $n_{\fctpsi(i)} \geq 3^{-\borne} n_i$, and therefore 
$\limsup n_{i+1}/n_i \leq 2 - 3^{-\borne}$.

Now, assume $\limsup \frac{n_{i+1}}{n_i} < 2 - \eps$ with $\eps > 0$.
Then obviously we have $\fctpsi(i) \in \Netoile$ for $i$ sufficiently large. Moreover
Lemma \ref{lemdelemtri} yields, for $k$ large enough:
\begineq \label{eqintermedpropbo}
n_{t_k - 1} < (1 - \eps) n_{t_k}. 
\eneq
Let $s_k = t_k - t_{k-1}$; the previous inequality yields
$$n_{t_k} = n_{t_{k-1}} + s_k (n_{t_k} - n_{t_k-1}) > n_{t_{k-1}} + s_k \eps n_{t_k},  $$ 
therefore $s_k \leq 1/\eps$: the sequence $(s_k)$ is bounded.
Now for $i$ large enough, and $\lambda$ such that
$\fctpsi(i) < t_{k-\lambda} < t_k \leq i$, Equation \eqref{eqintermedpropbo} gives
$$\eps n_i \leq n_{i+1} - 2 n_i = n_{\fctpsi(i)} \leq (1 - \eps)^{\lambda + 1} n_i$$
that is, an upper bound on $\lambda$. This concludes the proof of the lemma.\qed

\subsection{An Independence Property}  \label{subsecpreuvepropindep}

This Section is devoted to a proof of 
Proposition \ref{propindep}, which  means that $\densi(\fctpsi)$ 
can be defined using 
 the asymptotic behavior of any solution $(n'_i)$
of the associated recurrence relation: the initial values 
of the sequence $(n_i)$ do not matter.

\bigskip

Let $(n'_i)$ be as in the statement of Proposition \ref{propindep} 
(see \S \ref{subsecdemtheoqcq}).
Lemma \ref{lembornee} provides an integer
 $\borne$ such that $\fctpsi(i) \in \Netoile$ and 
$\fctpsi(i) \geq i - \borne$ for all sufficiently large $i$. 
Let $\difn'_i = n'_{i+1} - n'_i$. The recurrence relation 
satisfied by $(n'_i)$ 
 yields (for $i$ sufficiently large)
$$\difn'_{i+1} = \sum_{j = \fctpsi(i+1)} ^{i} \difn'_j . $$
This implies that $(\difn'_i)$ is non-decreasing for $i$ sufficiently large, 
and tends to infinity as $i$ tends to infinity
(except in the special case where $\fctpsi(i) = i-1$ for any sufficiently
large $i$, which is easily dealt with). Moreover, the same properties
hold for the sequence $\difn_i = n_{i+1}-n_i$.
 Now, let
$$\quo_i = \difn_i / \difn'_i. $$
The relation above yields, for $i$ sufficiently large:
$$\sum_{j= \fctpsi(i+1)} ^{i} \quo_j \difn'_j =
\sum_{j= \fctpsi(i+1)} ^{i} \difn_j = \difn_{i+1} 
= \quo_{i+1} \difn'_{i+1} = 
 \quo_{i+1} \sum_{j= \fctpsi(i+1)} ^{i} \difn'_{j} , $$
hence
\begineq \label{eqencadrtquo}
\quo_{i+1} = \sum_{j= \fctpsi(i+1)} ^{i} \frac{\difn'_j}{\difn'_{\fctpsi(i+1)} + \ldots + 
\difn'_i} \quo_j . 
\eneq
Now let $\imin$ be sufficiently large, and for $i \geq \imin$ let
$I_i$ be the convex hull of $\quo_{i+1-\borne}$, 
\ldots, $\quo_i$ (that is, the smallest segment in $\Rplusetoile$ that
contains these points). We shall deduce from \eqref{eqencadrtquo}
the following claim (where $\lgr{I}$ denotes the length of a segment $I$):
$$I_{i+\borne} \inclus I_i \mbox{ and }
\lgr{I_{i+\borne}} \leq (1 - \borne^{-\borne}) \lgr{I_i} \mbox{ for any }
i \geq \imin. $$
If the claim holds then the intersection of all $I_i$'s is reduced to a positive real number, 
which is the limit of the sequence $(\quo_i)$. As 
$\difn'_i$ tends to infinity with $i$, it is a classical consequence that
$$\frac{n_i}{n'_i} = \frac{n_0 + \difn_0 + \difn_1 + \ldots +
 \difn_{i-1}}{n'_0 + 
\difn'_0 + \difn'_1 + \ldots + \difn'_{i-1}}$$
converges to the same limit. This concludes the proof of the proposition --
if the claim holds.

To prove the claim, 
write $I_i = [\quo_i - \distamoins_i, \quo_i + \distaplus_i]$ and
notice that the right handside of \eqref{eqencadrtquo}
is a linear combination of $\quo_{i+1-\borne}$, 
\ldots, $\quo_i$ with non-negative coefficients. Moreover, the coefficient
of $\quo_i$ is at least $1/\borne$. Therefore bounding 
 $\quo_j$ in \eqref{eqencadrtquo} from below by
$\quo_i - \distamoins_i$ 
(resp. from above by $\quo_i + \distaplus_i$)
for $j \neq i$ shows that
$$\quo_{i+1} \in [\quo_i - (1 - 1/\borne) \distamoins_i, 
\quo_i + (1 - 1/\borne) \distaplus_i]. $$
Applying this result inductively yields, for any $\ell \geq 0$:
$$\quo_{i+\ell} \in [\quo_i - (1 - \borne^{-\ell}) \distamoins_i, \quo_i + 
(1 - \borne^{-\ell}) \distaplus_i] \inclus I. $$
This proves the claim, thereby concluding the proof of  Proposition \ref{propindep}.\qed

\subsection{A Special Property of $\sqrt3$}  \label{subsecpreuvelemdecalcul}

In this Section, we prove  Lemma \ref{lemdecalcul} stated in \S \ref{subsecpreuvetheosqrttrois}.

\bigskip

Let $\fctpsi$ be a reduced function such that $\densi(\fctpsi) < \sqrt3$, and 
  $(n_i)$ be the  associated sequence,   as in \S \ref{subsecappendicedef}. For any $k \geq 1$ and any $i \in \{t_k, \ldots, t_{k+1}-2\}$, we have $\fctpsi(i+1)=i$ hence
$n_{i+2}-n_{i+1} =   n_{i+1} - n_i$. Let us denote by $\qaa_k$ the common value of
$n_{i+1} - n_i$ for $i \in \{t_k, \ldots, t_{k+1}-1\}$. 
Excluding the case where the family $(t_k)$ is finite, it is clear that
the sequence $(\qaa_k)$ is increasing.
Let $\delta < \sqrt3$ and $K \geq 2$ be such that 
$n_{i+1} \leq \delta n_i$ 
(hence $\fctpsi(i) \in \Netoile$)
for all $i \geq t_K$.

Assume there is an index $k > K$ such that $\fctpsi(t_k) < t_{k-1} - 1$.

First of all, let us write (using Lemma \ref{lemdelemtri} with $i=t_{k-1}$)
\begineq \label{eqintermedunpropcalc}
n_{\fctpsi(t_k)} = 2 n_{t_k} - n_{t_k+1} \geq (2 - \delta) n_{t_k} 
\geq (2-\delta) (n_{t_{k-1}}+n_{t_{k-1}-1}). 
\eneq
Substracting $\qaa_{k-2} = n_{t_{k-1}} - n_{t_{k-1}-1}$ from
 this inequality yields (since $\delta \geq 1$):
\begineq \label{eqintermedpropcalc}
n_{\fctpsi(t_k)} - \qaa_{k-2} \geq (1-\delta) n_{t_{k-1}} + (3-\delta) n_{t_{k-1}-1}
\geq (3 - \delta^2) n_{t_{k-1}-1}. 
\eneq
Since $\delta < \sqrt3$, the right handside of \eqref{eqintermedpropcalc}
is positive. Therefore Lemma \ref{lemdelemtri} yields
$n_{\fctpsi(t_k)} > \qaa_{k-2} > n_{t_{k-2}-1}$ hence $\fctpsi(t_k) \geq t_{k-2}$. 
This inequality, together with the assumption $\fctpsi(t_k) < t_{k-1} - 1$, yields
$n_{t_{k-1}-1} - \qaa_{k-2} = n_{t_{k-1}-2} \geq n_{\fctpsi(t_k)}$. Combining this with
\eqref{eqintermedunpropcalc} gives
\begineq \label{eqntkak}
(2 \delta - 3) n_{t_{k-1}-1} \geq (3 - \delta) \qaa_{k-2}, 
\eneq
which implies, in particular, $\delta \geq 3/2$.
Now we also have 
$$\delta n_{t_{k-1}} \geq n_{t_{k-1}+1} \geq n_{t_{k-1}} + n_{t_{k-1}-1} = 2 n_{t_{k-1}} - \qaa_{k-2}, $$
hence 
$$\qaa_{k-2} \geq (2 - \delta) n_{t_{k-1}} \geq (2 - \delta) (\qaa_{k-2} + n_{t_{k-1}-1}). $$
Combining this relation with \eqref{eqntkak} yields $\delta \geq \sqrt3$, that is a contradiction.\qed

\newcommand{\url}{\texttt}
\providecommand{\bysame}{\leavevmode ---\ }
\providecommand{\og}{``}
\providecommand{\fg}{''}
\providecommand{\smfandname}{\&}
\providecommand{\smfedsname}{\'eds.}
\providecommand{\smfedname}{\'ed.}
\providecommand{\smfmastersthesisname}{M\'emoire}
\providecommand{\smfphdthesisname}{Th\`ese}

\bigskip

\bigskip

\hspace{-\parindent}St\'ephane Fischler

\hspace{-\parindent}\'Equipe d'Arithm\'etique et de G\'eom\'etrie Alg\'ebrique

\hspace{-\parindent}B\^atiment 425

\hspace{-\parindent}Universit\'e Paris-Sud

\hspace{-\parindent}91405 Orsay Cedex, France

\hspace{-\parindent}stephane.fischler@\null{}math.u-psud.fr


\begin{thebibliography}{10}

\bibitem{ABCD}
{\scshape J.~Allouche, M.~Baake, J.~Cassaigne {\normalfont \smfandname}
  D.~Damanik} -- {\og Palindrome complexity\fg}, \emph{Theoret. Comput. Sci.}
  \textbf{292} (2003), no.~1, p.~9--31.

\bibitem{BL}
{\scshape Y.~Bugeaud {\normalfont \smfandname} M.~Laurent} -- {\og Exponents of
  diophantine approximation and {S}turmian continued fractions\fg}, \emph{Ann.
  Inst. Fourier (Grenoble)} \textbf{55} (2005), p.~773--804.

\bibitem{Cassaigne}
{\scshape J.~Cassaigne} -- {\og Limit values of the recurrence quotient of
  {S}turmian sequences\fg}, \emph{Theoret. Comput. Sci.} \textbf{218} (1999),
  p.~3--12.

\bibitem{DamanikGhezRaymond}
{\scshape D.~Damanik, J.-M. Ghez {\normalfont \smfandname} L.~Raymond} -- {\og
  A palindromic half-line criterion for absence of eigenvalues and applications
  to substitution hamiltonians\fg}, \emph{Ann. H. Poincar\'e} \textbf{2}
  (2001), p.~927--939.

\bibitem{DS}
{\scshape H.~Davenport {\normalfont \smfandname} W.M.~Schmidt} -- {\og
  Approximation to real numbers by algebraic integers\fg}, \emph{Acta Arith.}
  \textbf{15} (1969), p.~393--416.

\bibitem{DeLuca}
{\scshape A.~{de Luca}} -- {\og Sturmian words: structure, combinatorics, and
  their arithmetics\fg}, \emph{Theoret. Comput. Sci.} \textbf{183} (1997),
  p.~45--82.

\bibitem{LucaMignosi}
{\scshape A.~{de Luca} {\normalfont \smfandname} F.~Mignosi} -- {\og Some
  combinatorial properties of {S}turmian words\fg}, \emph{Theoret. Comput.
  Sci.} \textbf{136} (1994), p.~361--385.

\bibitem{DJP}
{\scshape X.~Droubay, J.~Justin {\normalfont \smfandname} G.~Pirillo} -- {\og
  Episturmian words and some constructions of de {L}uca and {R}auzy\fg},
  \emph{Theoret. Comput. Sci.} \textbf{255} (2001), p.~539--553.

\bibitem{SFOttdio}
{\scshape S.~Fischler} -- {\og Palindromic prefixes and diophantine
  approximation\fg}, preprint arxiv math.NT/0509508.

\bibitem{CRASasim}
\bysame , {\og Spectres pour l'approximation d'un nombre r\'eel et de son
  carr\'e\fg}, \emph{C. R. Acad. Sci. Paris, Ser. I} \textbf{339} (2004),
  no.~10, p.~679--682.

\bibitem{HofKnillSimon}
{\scshape A.~Hof, O.~Knill {\normalfont \smfandname} B.~Simon} -- {\og Singular
  continuous spectrum for palindromic {S}chr\"odinger operators\fg},
  \emph{Commun. Math. Phys.} \textbf{174} (1995), p.~149--159.

\bibitem{JP}
{\scshape J.~Justin {\normalfont \smfandname} G.~Pirillo} -- {\og Episturmian
  words and episturmian morphisms\fg}, \emph{Theoret. Comput. Sci.}
  \textbf{276} (2002), p.~281--313.

\bibitem{Lothaire}
{\scshape M.~Lothaire} -- \emph{Algebraic combinatorics on words}, Encyclopedia
  of Mathematics and its Applications, no.~90, Cambridge University Press,
  2002.

\bibitem{MorseHedlund}
{\scshape M.~Morse {\normalfont \smfandname} G.A.~Hedlund} -- {\og Symbolic
  dynamics {II}: {S}turmian trajectories\fg}, \emph{Amer. J. Math.} \textbf{62}
  (1940), p.~1--42.

\bibitem{RauzySMF}
{\scshape G.~Rauzy} -- {\og Nombres alg\'ebriques et substitutions\fg},
  \emph{Bull. Soc. Math. France} \textbf{110} (1982), p.~147--178.

\bibitem{RauzyBordeaux}
\bysame , {\og Suites \`a termes dans un alphabet fini\fg}, \emph{S\'em.
  Th\'eor. Nombres Bordeaux} (1982-1983), p.~25.01--25.16.

\bibitem{RoyCRAS}
{\scshape D.~Roy} -- {\og Approximation simultan\'ee d'un nombre et de son
  carr\'e\fg}, \emph{C. R. Acad. Sci. Paris, Ser. I} \textbf{336} (2003),
  p.~1--6.

\end{thebibliography}
\end{document}